\pdfoutput=1
\RequirePackage{ifpdf}
\ifpdf 
\documentclass[pdftex]{sigma}
\else
\documentclass{sigma}
\fi

\newtheorem{Theorem}{Theorem}[section]
\newtheorem{Corollary}[Theorem]{Corollary}
\newtheorem{Lemma}[Theorem]{Lemma}

\newtheorem{Conjecture}[Theorem]{Conjecture}
\newtheorem{Question}[Theorem]{Question}

{\theoremstyle{definition}
\newtheorem{Definition}[Theorem]{Definition}

\newtheorem{Example}[Theorem]{Example}
\newtheorem{Remark}[Theorem]{Remark}
\newtheorem{Notation}[Theorem]{Notation}
}

\newcommand{\source}[1]{\operatorname{reactant}(#1)}
\newcommand{\target}[1]{\operatorname{product}(#1)} 
\newcommand{\flux}[1]{\operatorname{f\/lux}(#1)} 

\newcommand{\init}[2]{\operatorname{init}_{#1}\left(#2\right)}
\newcommand{\supp}[2]{\operatorname{supp}_{#1}\left(#2\right)}

\newcommand\ol[1]{{\overline {#1}}}
\renewcommand\iff{\Leftrightarrow}

\begin{document}
\allowdisplaybreaks

\renewcommand{\PaperNumber}{025}

\FirstPageHeading

\ShortArticleName{A Projection Argument for Dif\/ferential Inclusions}

\ArticleName{A Projection Argument for Dif\/ferential Inclusions,\\ with Applications to Persistence of Mass-Action\\
Kinetics}

\Author{Manoj GOPALKRISHNAN~$^\dag$, Ezra MILLER~$^\ddag$ and Anne SHIU~$^\S$}

\AuthorNameForHeading{M.~Gopalkrishnan, E.~Miller and A.~Shiu}

\Address{$^\dag$~School of Technology and Computer Science, Tata Institute of Fundamental Research,
\\
\hphantom{$^\dag$}~1~Homi Bhabha Road, Mumbai 400 005, India}
\EmailD{\href{mailto:manojg@tifr.res.in}{manojg@tifr.res.in}} \URLaddressD{\url{http://www.tcs.tifr.res.in/~manoj/}}

\Address{$^\ddag$~Department of Mathematics, Duke University, Box 90320, Durham, NC 27708-0320, USA}
\EmailD{\href{mailto:ezra@math.duke.edu}{ezra@math.duke.edu}} \URLaddressD{\url{http://www.math.duke.edu/~ezra/}}

\Address{$^\S$~Department of Mathematics, University of Chicago,
\\
\hphantom{$^\S$}~5734 S.~University Avenue, Chicago, IL 60637, USA}
\EmailD{\href{mailto:annejls@math.uchicago.edu}{annejls@math.uchicago.edu}}
\URLaddressD{\url{http://math.uchicago.edu/~annejls/}}

\ArticleDates{Received August 07, 2012, in f\/inal form March 23, 2013; Published online March 26, 2013}

\Abstract{Motivated by questions in mass-action kinetics, we introduce the notion of \emph{vertexical family} of
dif\/ferential inclusions.
Def\/ined on open hypercubes, these families are characterized by particular good behavior under projection maps.
The motivating examples are certain families of reaction networks~-- including reversible, weakly reversible,
endotactic, and \emph{strongly endotactic} reaction networks~-- that give rise to vertexical families of mass-action
dif\/ferential inclusions.
We prove that vertexical families are amenable to structural induction.
Consequently, a~trajectory of a~vertexical family approaches the boundary if and only if either the trajectory
approaches a~vertex of the hypercube, or a~trajectory in a~lower-dimensional member of the family approaches the
boundary.
With this technology, we make progress on the global attractor conjecture, a~central open problem concerning
mass-action kinetics systems.
Additionally, we phrase mass-action kinetics as a~functor on reaction networks with variable rates.}

\Keywords{dif\/ferential inclusion; mass-action kinetics; reaction network; persistence; global attractor conjecture}

\Classification{34A60; 80A30; 92C45; 37B25; 34D23; 37C10; 37C15; 92E20; 92C42; 54B30; 18B30 }

\section{Introduction}
\label{sec:intro}

\looseness=-1
The global attractor conjecture has been a~central open problem in reaction network theory since its formulation by
Horn in 1974~\cite{horn74dynamics}.
It asserts that any complex-balanced mass-action kinetics system of ordinary dif\/ferential equations with positive
initial conditions possesses a~globally attracting stationary point in each stoichiometric compatibility class (see
Conjecture~\ref{c:GAC} for a~more precise statement).
It is well-known that this conjecture is implied by Feinberg's persistence conjecture~\cite[Remark~6.1.E]{Fein87},
a~version of which asserts the following: for weakly reversible networks taken with mass-action kinetics, no species
asymptotically becomes extinct or unbounded.

An a~priori special case of the persistence conjecture asserts that at least one species survives asymptotically.
In this paper, we show that this a~priori special case in fact implies the entire persistence conjecture
(Corollary~\ref{cor:repelled}).
As a~consequence, if the persistence conjecture is false, then in every minimal counterexample each species becomes
either extinct or unbounded.
It follows that the persistence conjecture in dimension $n$ implies the global attractor conjecture in dimension $n+1$
(Theorems~\ref{t:n=>n+1} and~\ref{t:bdd_tra_orig_rep=>conjs}).

Pantea~\cite{Pantea}, building on earlier work by Craciun, Nazarov, and Pantea~\cite{CNP}, used the persistence
conjecture in two dimensions to prove the global attractor conjecture in three dimensions.
Pantea's work relied in part on projecting trajectories to lower-dimensional faces.
Here, we generalize this projection argument in two ways.
First, our main result, Theorem~\ref{t:vertexical}, applies in arbitrary dimensions.
Second, our results hold not only for mass-action kinetics networks, but for certain families of dif\/ferential
inclusions that we call \emph{vertexical} (Def\/inition~\ref{d:vertexical}).
The notion of vertexical family makes precise the essential structure required of a~family of dynamical systems on
hypercubes of varying dimensions to permit a~structural induction argument of this sort.
Theorem~\ref{t:vertexical} shows that a~trajectory in a~vertexical family approaches the boundary if and only if either
it approaches a~vertex of the hypercube, or a~lower-dimensional trajectory in the family approaches the boundary.

Vertexical families of dif\/ferential inclusions arise naturally in reaction network theory by way of mass-action
kinetics or, more generally, power-law dynamics that are considered in biochemical systems theory
(Remark~\ref{rmk:genl_crn}), on networks that are reversible, weakly reversible, endotactic, strongly endotactic
(Def\/inition~\ref{def:endotactic}.\ref{def:strong_endotactic}), and so on.
We prove that these networks, and more generally, \emph{projective classes} of networks
(Def\/inition~\ref{def:projective}), give rise to vertexical families of mass-action dif\/ferential inclusions
(Theorem~\ref{t:proj-vertexical} and Corollary~\ref{cor:isvertexical}).

In the course of proving this result, we are led to view mass-action kinetics as a~functor (Theo\-rem~\ref{thm:functor}).
No more category theory is required in this paper beyond the def\/inition of a~functor.
Functoriality itself is used as a~convenient shorthand for a~list of properties spelled out at the beginning of
Section~\ref{sec:funct}.
The use of this shorthand clarif\/ies the concept of vertexical family and suggests that other questions concerning
mass-action kinetics systems may be amenable to structural induction~(Question~\ref{qopen}).
Section~\ref{sec:implications} discusses the implications of our results for persistence of mass-action kinetics
systems.

\section{Dynamical properties of dif\/ferential inclusions}\label{s:dynamical}

In this section, we recall certain dynamical properties of dif\/ferential inclusions def\/ined on ma\-ni\-folds.
For background on manifolds, see~\cite{Lee02smooth}.
All manifolds considered here have f\/inite dimension.
For background on dif\/ferential inclusions, see~\cite{Aubin}.
\begin{Definition}\label{def:diff_incl}
Let $M$ be a~smooth manifold with tangent bundle $\pi_M:TM\to M$.
A~\emph{dif\-fe\-ren\-tial inclusion} on $M$ is a~subset $X\subseteq TM$.
\end{Definition}

\begin{Example}\label{e:vectorfield}
The simplest dif\/ferential inclusions on~$M$ are vector f\/ields on~$M$.
The subset $X\subseteq TM$ for a~given vector f\/ield is the image of the corresponding section $M\hookrightarrow TM$.
\end{Example}

\begin{Definition}\label{def:trajectory}
Fix a~dif\/ferential inclusion $X$ on a~smooth manifold~$M$.
Let $I\subseteq\mathbb{R}_{\geq 0}$ be a~nonempty interval (in particular, connected) containing its left endpoint.
A~dif\/fe\-ren\-tiab\-le curve $f:I\to M$ is a~\emph{trajectory} of $X$ if the tangent vectors to the curve lie in~$X$.
An unbounded interval is a~\emph{ray}.
A trajectory def\/ined on a~ray \emph{eventually} has a~property $P$ if there exists $T>0$ such that property $P$ holds
for the function for all $t\geq T$.
\end{Definition}

Let $\ol M$ be a~smooth manifold with corners.
That is, $\ol M$ is a~space locally modeled on the closed nonnegative orthant~\cite[p.~363]{Lee02smooth}.
Then $\partial\ol M$ denotes the boundary of~$\ol M$, which is the set of points of $\ol M$ that are not in the
relative interior of~$\ol M$.
The relative interior $\ol M\setminus\partial\ol M$ of $\ol M$
is a~smooth manifold~\cite[p.~386, Examples~14--19]{Lee02smooth}.

\begin{Definition}
\label{d:persistent}
Let $\ol{M}$ be a~smooth manifold with corners with relative interior $M=\ol M\setminus\partial\ol M$, and let
$V\subseteq\partial\ol M$ be a~subset of the boundary.
\begin{enumerate}\itemsep=0pt
\item A dif\/ferential inclusion $X\subseteq TM$ is \emph{persistent relative to~$V$} if the closure in~$\ol M$ of
every trajectory of~$X$ is disjoint from the closure~$\ol V$ of $V$ in $\ol M$.

\item A dif\/ferential inclusion $X\subseteq TM$ is \emph{repelled by $V$} if for every open set $O_1\subseteq\ol M$
with $\ol V\subseteq O_1$, there exists a~smaller open set $O_2\subseteq O_1$ with $\ol V\subseteq O_2$ such that for
every trajectory $f:I\to M$ of~$X$, if $f(\inf{I})\notin O_1$ then $f(I)\cap O_2$ is empty; in other words, if the
trajectory begins outside of~$O_1$, then the trajectory never enters $O_2$.

\item If $\ol M$ is compact, then a~dif\/ferential inclusion $X\subseteq TM$ is \emph{permanent} if it is persistent
and there is a~compact subset $\Omega\subseteq M$ such that for every ray~$I$, every trajectory of~$X$ def\/ined on~$I$
is eventually contained in~$\Omega$.
\end{enumerate}
More generally, a~set~$\mathcal X$ of dif\/ferential inclusions on~$M$ is \emph{persistent relative to~$V$},
\emph{repelled by~$V$}, or \emph{permanent} if every member $X\in\mathcal X$ has the corresponding property.
\end{Definition}
\begin{Definition}\label{d:persistent'}
Dif\/ferential inclusions that are persistent relative to the boundary $\partial\ol M$ are simply called
\emph{persistent} and similarly for \emph{repelled}.
A collection of dif\/ferential inclusions, possibly on a~family of dif\/ferent manifolds with corners, is
\emph{persistent}, \emph{permanent}, or \emph{repelled} if each dif\/ferential inclusion in the collection has the
respective property.
\end{Definition}
\begin{Remark}\label{r:repelled==>persistent}
Any dif\/ferential inclusion repelled by~$V$ is also persistent relative to~$V$.
The converse is false in general because dif\/ferent trajectories starting outside $O_1$ could get arbitrarily close
to~$V$.
However, certain extra conditions could guarantee that the original dif\/ferential inclusion is repelled by~$V$.
For example, suppose~$\ol M$ is compact and that the dif\/ferential inclusion~$X$ has a~continuous extension~$\ol X$
to~$\ol M$.
Assume further that every trajectory of~$\ol X$ starting in~$\partial\ol M$ but outside~$\ol V$ has its closure
disjoint from~$\ol V$.
If the projection $\ol X\to\ol M$ is suf\/f\/iciently nice~-- we are unsure what conditions to impose, but we have in
mind properness~-- then it should be possible to conclude that~$X$ is repelled~by~$V$.
\end{Remark}

\begin{Remark}
\label{r:permanent}
A dif\/ferential inclusion that is permanent need not be repelled by~the boundary of~$\ol M$.
The reason appeared already in Remark~\ref{r:repelled==>persistent}: dif\/ferent trajectories starting outside~$O_1$
could get arbitrarily close to the boundary (on the way to ending up in~$\Omega$).
Conversely, a~dif\/ferential inclusion repelled by the boundary need not be permanent, even if $\ol M$ is compact,
because $O_2$ might necessarily be smaller when $O_1$ is smaller.
That is, trajectories that start closer to the boundary could eventually remain closer to the boundary; see
Example~\ref{e:repelled}.
\end{Remark}
\begin{Example}
\label{e:repelled}
Fix a~dif\/ferential inclusion~$X$ on a~planar disk~$\ol M$ whose trajectories form concentric circles about its center.
$X$ is repelled by the boundary circle $V=\partial\ol M$.
Indeed, if $O_1\subseteq\ol M$ is an open set containing~$V$, then the compact set $\ol M\setminus O_1$ achieves
a~maximum radius~$r$ from the center, so we can take $O_2$ to be the set of all points in~$\ol M$ of radius $>r$.
\end{Example}
\begin{Lemma}
\label{l:normal}
If the differential inclusion $X$ in Definition~{\rm \ref{d:persistent}} is persistent, then for every trajec\-tory $f:I\to
M$ of~$X$, there exist disjoint open sets~$O_f$ and~$O_V$ in~$\ol M$ containing the closures in~$\ol M$ of~$f(I)$
and~$V$ respectively.
\end{Lemma}

\begin{proof}
A manifold with corners is metrizable, and hence it is a~normal Hausdorf\/f space.
\end{proof}

\begin{Remark}
\label{rmk:def_pers}
Consider a~dif\/ferential inclusion def\/ined on a~positive orthant $M=\mathbb{R}_{>0}^n$.
Distinct partial compactif\/ications $\ol M$ of $M$ yield distinct notions of persistence.
In the case of $\ol M=\mathbb{R}_{\geq 0}^n$, even with Lemma~\ref{l:normal}, our def\/inition of persistence is weaker
than the standard def\/inition~\cite{Freed} that requires each coordinate of a~trajectory to remain bounded away
from~$0$ for all time.
However, the mass-action dif\/ferential inclusions that we consider are viewed in the compactif\/ication $\ol
M=[0,\infty]^S$ (Remark~\ref{rmk:m-a_in_cpc}), so in the context of reaction network theory, our def\/inition of
persistence is in fact stronger than the standard def\/inition because coordinates of trajectories not only must be
bounded away from zero but also must avoid going to inf\/inity.

\looseness=-1
Mathematically, the def\/inition we adopt has the advantage of being purely topological, so it behaves well under
homeomorphism.
Our def\/inition also allows, if required, to separate the usual concept of persistence into the two questions of
whether trajectories are bounded and whether $\omega$-limit points exist on the boundary.
Both properties are conjectured to hold for weakly reversible and, more generally, endotactic (see
Def\/inition~\ref{def:endotactic}) reaction networks.
When trajectories are bounded, our def\/inition of persistence is equivalent to the standard one.
Indeed, the mass-action dif\/ferential inclusions that we introduce later are viewed in the compactif\/ication
$[0,\infty]^{S}$, so our def\/inition of persistence automatically implies boundedness of trajectories; see
Remark~\ref{rmk:m-a_in_cpc}.
\end{Remark}

\begin{Remark}
\label{rmk:repel_dist}
Suppose a~manifold $M$ has a~metric $\mathrm{d}$, and $V$ is a~compact subset of~$M$.
A~dif\/ferential inclusion~$X$ is repelled by~$V$ if for every $d_1\in\mathbb{R}_{>0}$, there exists
$d_2\in\mathbb{R}_{>0}$ such that every trajectory $f:I\to M$ of~$X$ starting with $\mathrm{d}(f(\inf{I}),V)\geq d_1$
maintains $\mathrm{d}(f(I),V)\geq d_2$.
\end{Remark}

\begin{Remark}
\label{r:repelled}
Our notion of ``repelled'' is new, motivated by the requirements of Theorem~\ref{t:vertexical}.
The motivations are further explained in Remark~\ref{r:why_can't_replace_repelling_with_persistent}.
Cognate but distinct concepts bearing similar names have been def\/ined by others.
Anderson and Shiu def\/ine a~boundary face to have a~``repelling neighborhood'' if there is a~neighborhood of the face
such that whenever a~trajectory enters that neighborhood, it can get no closer to the face while remaining in that
neighborhood~\cite{AndersonShiu10}.
Banaji and Mierczynski def\/ine a~``repelling face'' as a~boundary face for which any trajectory that begins in that
face immediately exits the face into the interior of the relevant invariant set~\cite{BanajiMier}.
Neither of these concepts is adequate for our purposes.
\end{Remark}
\begin{Remark}
\label{r:subsetTrajectory}
Suppose a~dif\/ferential inclusion~$X$ is a~subset of a~persistent dif\/ferential inclusion~$Y$.
Then $X$ must be persistent, since each of its trajectories is a~trajectory of~$Y$.
More generally, consider properties~$P$ of dif\/ferential inclusions~$X$ of the form ``$P(f)$ holds for all
trajectories~$f$ of~$X$.'' If property~$P$ is true for a~dif\/ferential inclusion~$Y$, and a~dif\/ferential
inclusion~$Z$ factors through~$Y$, then property $P$ is true for~$Z$ as well.
\end{Remark}

\section{Vertexical families of dif\/ferential inclusions}
\label{s:vertexical}

A mass-action kinetics system is naturally def\/ined on the positive orthant $\mathbb{R}_{>0}^S$ corresponding to the
space of concentrations of species.
In this section we work instead with open hypercubes $(0,1)^S$, and not directly with positive orthants themselves.

\looseness=-1
To justify this choice, f\/irst we argue that nothing is lost by working with hypercubes.
Open hypercubes $(0,1)^S$ are dif\/feomorphic to positive orthants $\mathbb{R}_{>0}^S$, so dif\/ferential inclusions can
be transferred from one space to the other by f\/ixing a~dif\/feomorphism and using its Jacobian (see
Section~\ref{subsec:functor}).
Additionally, properties such as persistence are def\/ined topologically on the tangent bundle~-- and therefore
invariant under dif\/feomorphism~-- so they can be analyzed on either space.

An advantage of working with open hypercubes is that they have natural ``cubical'' compactif\/ications $[0,1]^S$ that
appear to be optimally relevant in the context of mass-action kinetics.
Alternatively, we could have achieved a~cubical compactif\/ication by considering the hypercubes $[0,\infty]^S$.
Nevertheless, there is a~stylistic advantage to working with the hypercubes $[0,1]^S$: we can treat symmetrically the
cases where a~species concentration goes to inf\/inity or to zero.
This makes some of our def\/initions more transparent, and the structural induction becomes cleaner.
\begin{Definition}
\label{def:diff_incl_fam}
For any f\/inite nonempty set~$S$, let $\mathcal D_S$ be the set of all dif\/ferential inclusions on the open hypercube
$(0,1)^S$.
Fix a~collection $\mathcal S$ of f\/inite nonempty sets.
If $\mathcal X_S$ is a~set of dif\/ferential inclusions on $(0,1)^S$ for each $S\in\mathcal S$, then the collection
$\mathcal X=\{\mathcal X_S\subseteq\mathcal D_S\}_{S\in\mathcal S}$ of sets~$\mathcal X_S$ is a~\emph{family of
differential inclusions on open hypercubes indexed by~$\mathcal S$}.
\end{Definition}

\subsection{Definitions concerning hypercubes}
\label{sub:hypercubes}

The def\/inition of vertexical families requires some preliminary notation on hypercubes.
\begin{Notation}
Let $S$ be a~f\/inite nonempty set.
\begin{enumerate}\itemsep=0pt
\item For every $i\in S$, let $e_i\in\mathbb{R}^S$ be the standard basis vector indexed by~$i$; that is,
$e_i:S\to\mathbb{R}$ sends $i$ to $1$ and $S\setminus\{i\}$ to~$0$.
\item For each subset $U\subseteq S$ and vertex $x\in\{0,1\}^S$ of the hypercube $[0,1]^S$, let
\begin{gather*}
F_U(x)=\big(x+\operatorname{span}\{e_i\,|\, i\in U\}\big)\cap[0,1]^S
\end{gather*}
denote the face of the hypercube $[0,1]^S$ along $U$ at vertex $x$.
\item For $p=\sum\limits_{i\in S}p_i e_i\in\mathbb{R}^S$, let $|p|=\sqrt{\sum\limits_{i\in S}p_i^2}$ denote its Euclidean norm.
\item For subsets $P,Q\subseteq\mathbb{R}^S$, denote by $\mathrm{d}(P,Q)=\inf\big\{|p-q|\,\big|\, p\in P~\text{and}~q\in
Q\big\}$ the distance between them.
\end{enumerate}
\end{Notation}

\begin{Remark}
Our notation for faces dif\/fers from that in related references~\cite{Anderson08, AndersonShiu10, angeli2007Petri,
TDS, Siegel}.
What those works call $F_U$ is close to what we call $F_{S\setminus U}({0})$, where ${0}$ denotes the origin.
This correspondence is not perfect; the sets $F_U$ in the related works are faces of the nonnegative orthant
$\mathbb{R}^S_{\geq0}$, whereas here the sets $F_U(x)$ denote faces of hypercubes.
\end{Remark}

In the context of reaction networks, $S$ indexes the set of reacting chemical species; a~chemical complex, being
a~linear combination of these species, is therefore viewed as a~vector in $\mathbb{R}^S$ (see
Def\/inition~\ref{def:crn}), which has preferred basis vectors~$e_i$ for $i\in S$.
\begin{Definition}
\label{d:shrink}
Fix a~f\/inite nonempty set~$S$.
For a~face $F$ of $[0,1]^S$ and a~real number $\eta>1/2$, the \emph{centered shrinking}
\begin{gather*}
\eta F=\big\{x\in F\,\big|\,\mathrm{d}(x,\partial F)\geq(1-\eta)/2\big\}
\end{gather*}
of~$F$ is the set of points in $F$ whose distance from the boundary $\partial F$ is at least $(1-\eta)/2$.
\end{Definition}

\begin{Example}
\label{e:shrinking}
A centered shrinking of the rightmost face of the $3$-cube looks as follows,
$$
\includegraphics[scale=0.85]{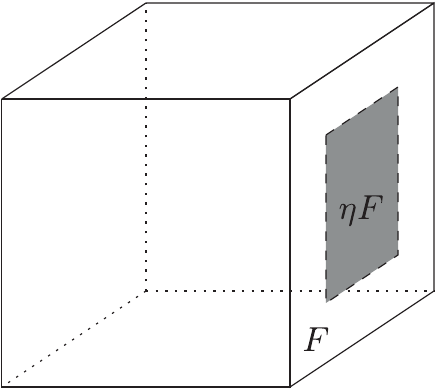}
$$
where the inner shaded square is centered in~$F$ and has side length $\eta$ times that of~$F$.
\end{Example}
\begin{Definition}
\label{d:pile}
Fix a~f\/inite nonempty set~$S$ and a~subset $U\subseteq S$.
Let $P\subseteq[0,1]^S$, and suppose that $\varepsilon\in(0,1/2)\subseteq\mathbb{R}$.
The \emph{$\varepsilon$-pile} of the subset $P$ along $U$ is the~set
\begin{gather*}
\operatorname{pile}(P,\varepsilon,U):=
\bigg\{x+\sum_{i\in U}\varepsilon_i e_i\,\Big|\,
x\in P~\text{and}~-\varepsilon\leq\varepsilon_i\leq\varepsilon~\text{for all}~i\in U\bigg\}\cap[0,1]^S.
\end{gather*}
\end{Definition}

\begin{Definition}
\label{d:thick}
Fix a~f\/inite nonempty set~$S$ and a~proper face $F$ of the hypercube $[0,1]^S$ containing a~vertex~$x$.
Let $U\subseteq S$ be such that $F=F_{S\setminus U}(x)$.
For a~real number $\varepsilon\in(0,1/2)$, the \emph{$\varepsilon$-block} $F_{\varepsilon}$ is
$\operatorname{pile}\big((1-2\varepsilon)F,\varepsilon,U\big)$, the $\varepsilon$-pile along $U$ of the centered
shrinking $(1-2\varepsilon)F$.
\end{Definition}

\begin{Example}
\label{e:block}
If $\varepsilon=(1-\eta)/2$ in Example~\ref{e:shrinking}, the block $F_\varepsilon$ looks like the following:
$$
\includegraphics[scale=0.85]{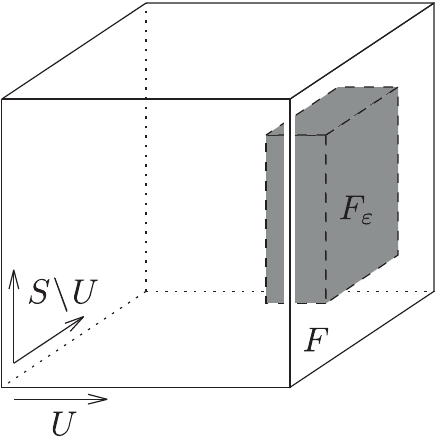}
$$
Note that the face~$F$ is orthogonal to the basis vector indexed by~$U$, and the thickness~$\varepsilon$ of the block
$F_\varepsilon$ equals its distance from the edges of~$F$.
The vertex~$x$ in Def\/inition~\ref{d:thick} could equally well be any of the four vertices of~$F$.
\end{Example}
\begin{Remark}
The block $F_\varepsilon$ is a~closed subset of $[0,1]^S$.
Such sets are closely related to sets that Pantea~\cite{Pantea} denoted by $K_\varepsilon$, which he used for the
purpose of projecting trajectories, as we too do in the current paper.
\end{Remark}
\begin{Notation}
\label{n:projectionmap}
Let $f:U\to S$ be a~map of f\/inite sets, and view $\mathbb{R}^S$ as functions $S\to\mathbb{R}$.
Denote by $\pi_f:\mathbb{R}^S\to\mathbb{R}^U$ the linear projection that sends $v\in\mathbb{R}^S$ to $v\circ
f\in\mathbb{R}^U$.
When $U\subseteq S$ and $f$ is inclusion, we write $\pi_U$ for the projection map instead of~$\pi_f$.
\end{Notation}
\begin{Remark}{\samepage
If $U\subseteq S$ is any subset, then
\begin{enumerate}\itemsep=0pt
\item[1)] projection is surjective on the open hypercube: $\pi_U\big((0,1)^S\big)=(0,1)^U$, and

\item[2)] $\pi_U\big((1-2\varepsilon)F_{S\setminus U}(x)\big)$ is a~vertex of the hypercube $[0,1]^U$.
\end{enumerate}}
\end{Remark}
\begin{Example}
\label{e:projection}
The projection $\pi_U$ in Example~\ref{e:block} collapses the cube to the horizontal~edge

\centerline{\includegraphics[scale=0.9]{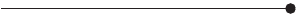}}

\noindent
that is $[0,1]=[0,1]^U$.
The projection takes $F$ as well as the subset $\eta F\subseteq F$ to the indicated right-hand vertex of the interval.
\end{Example}

\subsection{Definition of vertexical family and main result}
\label{sub:main}

The heuristic description of a~vertexical family of dif\/ferential inclusions begins by considering a~trajectory of
a~dif\/ferential inclusion in the family.
Suppose the trajectory remains near a~face of the hypercube.
The vertexical condition requires that, while the trajectory is near the face, the image of the trajectory under the
projection map collapsing that face be the trajectory of a~f\/ixed lower-dimensional dif\/ferential inclusion in the
family, up to time reparametrization.
We emphasize that only the part of the trajectory near the face and away from the boundary of the face is required to
be projectable.
\begin{Definition}\label{d:vertexical}
Let $\mathcal S$ be the set of all f\/inite nonempty subsets of the positive integers $\mathbb{Z}_{\geq1}$.
A family $\mathcal X=\{\mathcal X_S\}_{S\in\mathcal S}$ of dif\/ferential inclusions on open hypercubes indexed by
$\mathcal S$ is \emph{vertexical} if for each
\begin{itemize}\itemsep=-0.5pt
\item set $S\in\mathcal S$, \item dif\/ferential inclusion $X\subseteq T(0,1)^S$ in~$\mathcal X_S$, \item proper
nonempty subset $U\subseteq S$, and \item face $F=F_{S\setminus U}(x)$ of $[0,1]^S$,
\end{itemize}
there is $\varepsilon'>0$ such that for every $\varepsilon\in(0,\varepsilon')$, some dif\/ferential inclusion
$Y\in\mathcal X_U$ has the property that for every trajectory $f:I\to[0,1]^S$ of~$X$ with image in the
block~$F_\varepsilon$, there~exist
\begin{itemize}\itemsep=0pt
\item a~trajectory $g:J\to[0,1]^U$ of~$Y$, and \item an order-preserving continuous map $\alpha:I\to J$
\end{itemize}
such that $\pi_U\circ f=g\circ\alpha$.
\end{Definition}

Examples of vertexical families of dif\/ferential inclusions include those arising from reversible, weakly reversible,
endotactic, or strongly endotactic chemical reaction networks  (Def\/initions~\ref{def:network} and~\ref{d:functor});
this is the content of Corollary~\ref{cor:isvertexical}, the goal of Sections~\ref{sec:CRNT} and~\ref{sec:funct}.
Some nuances in the def\/inition are further discussed in Remark~\ref{rem:allotvxical}.

We now give a~def\/inition, followed by our main result on abstract vertexical families.
\begin{Definition}
\label{d:charged}
Fix a~f\/inite set $S$ and an index set $R\subseteq\mathbb{Z}_{\geq1}$, called the \emph{repulsing index set}.
Embed the hypercube $[0,1]^{R\cap S}$ into the hypercube $[0,1]^S$ as the face $[0,1]^{R\cap S}\times\{0\}^{S\setminus
R}$.
A~vertex of $[0,1]^S$ is \emph{charged} if it lies in $[0,1]^{R\cap S}$.
A~face $F$ of $[0,1]^S$ is \emph{opposite} if $F\cap[0,1]^{R\cap S}$ is empty.
The \emph{charged set} is the set $[0,1]^{R\cap S}\cap\partial[0,1]^S$.
\end{Definition}

In practice, $R$ and $S$ are both subsets of a~f\/ixed set, and $R$ need not be f\/inite.
If $R$ is empty, then by convention $[0,1]^{R\cap S}\times\{0\}^{S\setminus R}$ is the origin.
Thus the origin is always~charged.
The charged set equals $[0,1]^{R\cap S}$ unless $R\supseteq S$, in which case the charged set is $\partial[0,1]^S$.
\begin{Theorem}
\label{t:vertexical}
Fix a~vertexical family $\mathcal X=\{\mathcal X_S\}_{S\in\mathcal S}$ on open hypercubes indexed by the set $\mathcal
S$ of all finite nonempty subsets of the positive integers $\mathbb{Z}_{\geq1}$, and a~repulsing index set
$R\subseteq\mathbb{Z}_{\geq1}$.
Assume that for every set $S\in\mathcal S$, every differential inclusion $X\in\mathcal X_S$ is
\begin{itemize}\itemsep=-0.5pt
\item persistent relative to the union of all opposite faces of\/~$[0,1]^S$ and \item repelled by the charged vertices
of its hypercube $[0,1]^S$.
\end{itemize}
Then
\begin{enumerate}\itemsep=-0.5pt
\item[$1.$] Every such differential inclusion~$X$ is persistent relative to the entire boundary $\partial[0,1]^S$ and
repelled by the charged set $[0,1]^{R\cap S}\cap\partial[0,1]^S$ of its hypercube.
\item[$2.$] Fix $S\in\mathcal S$ and $X\in\mathcal X_S$.
If, in addition, $X$ is repelled by the union of all opposite faces of\/ $[0,1]^S$, then $X$ is repelled by the
boundary $\partial[0,1]^S$.
\end{enumerate}
\end{Theorem}

\begin{Remark}
\label{r:chargedSet}
The dif\/ferential inclusion in Theorem~\ref{t:vertexical}.1 is repelled either by the entire boundary of its hypercube
(if $R\supseteq S$) or by the proper face $[0,1]^{R\cap S}$.
\end{Remark}
\begin{proof}[Proof of Theorem~\ref{t:vertexical}] Fix $S\in\mathcal S$.
Let $X\in\mathcal X_S$.
We prove that for every proper, positive-dimensional face $F=F_{S\setminus U}(x)$ of~$[0,1]^S$, the following two
claims hold.
\begin{enumerate}\itemsep=-0.5pt
 \item[A.] If $X$ is persistent relative to the boundary $\partial F$ of $F$, then $X$ is persistent relative to~$F$.
\item[B.] \looseness=-1 If $F$ is not an opposite face and $X$ is repelled by the boundary $\partial F$ of $F$, then $X$ is repelled by
$F$.
\end{enumerate}
The remainder of this proof has two components: we f\/irst explain how Claims~A and~B imply parts~1 and~2 of the
theorem, and then we prove the two claims.

To start, consider part~1 of the theorem.
The dif\/ferential inclusion $X$ is persistent relative to all vertices of $[0,1]^S$, because each vertex is either
charged or an opposite face.
Using this fact as a~base case, Claim~A implies, by induction on the dimension of~$F$, that $X$ is persistent relative
to every proper face of $[0,1]^S$.
Since the hypercube $[0,1]^S$ has only f\/initely many faces, $X$ is therefore persistent relative to the entire
boundary.
For repulsion, all vertices of $[0,1]^S$ that lie in the charged set are charged vertices (see
Remark~\ref{r:chargedSet}), so $X$ is repelled by all such vertices by hypothesis.
Using this fact as a~base case, Claim~B implies, by induction on the dimension of faces~$F$ of~$[0,1]^S$ that are in
the charged set (and thus are not opposite faces), that $X$ is repelled by the charged set.
Hence, part~1 of the theorem holds for $X$, after we prove the two claims.

As for part~2, now assume that $X$ is repelled by the union of opposite faces of\/~$[0,1]^S$.
We need that $X$ is repelled by non-opposite faces as well.
Each vertex of $[0,1]^S$ is either charged or an opposite face, so each vertex repels~$X$.
Using this fact as a~base case, Claim~B implies, by induction on the dimension of proper non-opposite faces~$F$, that
$X$ is repelled by every non-opposite face.

It remains to prove Claims~A and~B for a~proper, positive-dimensional face $F=F_{S\setminus U}(x)$ of~$[0,1]^S$.
If $F$ is an opposite face, then Claim~A holds by hypothesis and Claim~B is vacuous.

Now assume that $F$ is not an opposite face.
Assume that~$X$ is persistent relative to the bounda\-ry~$\partial F$ of the face.
Let $f:I\to(0,1)^S$ be a~trajectory of~$X$, and let $d_1=\mathrm{d}\big(f(\inf{I}),F\big)$ be the distance to~$F$ from
the initial point of the trajectory.
The goal is to exhibit $\varepsilon>0$ so that the trajectory remains at distance greater than~$\varepsilon$ from~$F$;
that is, $\mathrm{d}\big(f(I),F\big)\geq\varepsilon$.
Claim~A then follows as a~consequence.

By hypothesis, $X$ is persistent relative to the boundary~$\partial F$, so there exists $d_2>0$ such that
$\mathrm{d}\big(f(I),\partial F\big)\geq d_2$.
Decreasing~$d_2$ if necessary, assume that $d_2\leq d_1$ and that $d_2/2\leq\varepsilon'$, where $\varepsilon'>0$ is
such that trajectories in the block~$F_{\varepsilon'}$ can be projected (Def\/inition~\ref{d:vertexical}).

Consider the block~$F_{d_2/2}$ of~$F$.
If the image of the trajectory $f(I)$ fails to intersect the block~$F_{d_2/2}$, then by def\/ining $\varepsilon=d_2/2$
it follows that $\mathrm{d}\big(f(I),F\big)>\varepsilon$, and we are done.
Therefore, we can and do assume that $I'\subseteq I$ is a~maximal nonempty subinterval such that $f(I')\subseteq
F_{d_2/2}$, and let $\iota=f(\inf{I'})$ denote the corresponding initial point, as in the following illustration.
$$
\includegraphics[scale=0.85]{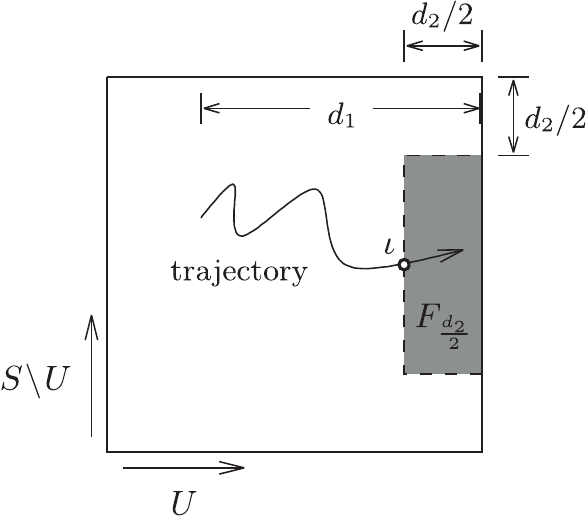}
$$

Note that $\mathrm{d}(\iota,F)=d_2/2$.
Indeed, by maximality of the interval~$I'$ the point $\iota$ lies on the boundary of~$F_{d_2/2}$, and the only boundary
face of $F_{d_2/2}$ that intersects the interior $(0,1)^S$ of the hypercube without also being contained in the
$d_2$-neighborhood of the boundary~$\partial F$ has constant distance $d_2/2$ from~$F$; see the following illustration.
$$
\includegraphics[scale=0.85]{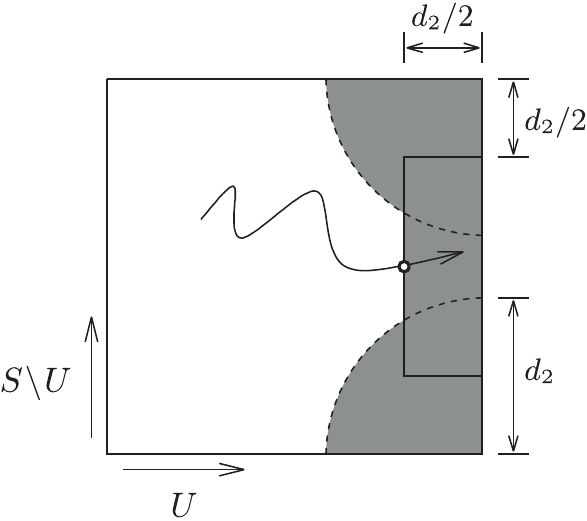}
$$

Because $\mathcal X$ is a~vertexical family, there exist a~dif\/ferential inclusion $Y$ in~$\mathcal X_U$, a~trajectory
$g:J\to[0,1]^U$ of~$Y$, and an order-preserving continuous map $\alpha:I'\to J$ such that $g\circ\alpha=\pi_U\circ f$
with domain~$I'$.
By def\/inition, $Y$ depends only on $X$ and~$F$, and not on the particular trajectory~$f$ or the subinterval $I'$.
To prove Claim~A for this face~$F$, it now suf\/f\/ices to show that there exists $\varepsilon>0$, depending only
on~$d_2$, such that $\mathrm{d}\big(g(J),\pi_U(F)\big)\geq\varepsilon$, for this claim implies that
$\mathrm{d}\big(f(I),F\big)\geq\varepsilon$, as desired.

Since $F$ is not an opposite face, by def\/inition it contains a~charged vertex, which we assume without loss of
generality is~$x$ (recall that $F=F_{S\setminus U}(x)$).
We claim that $y=\pi_U(F)$ is a~charged vertex of~$[0,1]^U$; that is, $y_i=0$ for all $i\in U\setminus R$.
Indeed, $j\in U$ implies $x_j=y_j$ because $x,y\in F=F_{S\setminus U}(x)$, and if additionally $j\notin R$, then
$x_j=0$ because $x$ is charged.

By hypothesis of the theorem, $Y$ is repelled by the charged vertex $y=\pi_U(F)$ of~$[0,1]^U$.
Hence there exists $\varepsilon>0$ such that all trajectories of $Y$ starting at distance $d_2/2$ away from the vertex
$\pi_U(F)$~-- and in particular, the trajectory~$g$, because $\iota$ has distance precisely $d_2/2$ from~$F$~-- never get
closer than $\varepsilon$ to the vertex~$\pi_U(F)$.
This choice of $\varepsilon$ depends only on the distance~$d_2/2$, not on any aspect of the particular trajectory~$f$;
its existence proves Claim~A.

\looseness=1
To prove Claim~B for this face~$F$, assume $X$ is repelled by~$\partial F$, and let $O_1\subseteq[0,1]^S$ be an open
set containing~$F$.
Using that $X$ is also persistent relative to~$\partial F$, repeat the argument above, but with $d_1>0$ now denoting
the distance from $[0,1]^S\setminus O_1$ to~$F$.
The assumption that~$X$ is repelled by~$\partial F$ implies that the value of~$d_2$ as found above depends only
on~$d_1$ and not on any aspect of any trajectory~$f$ that begins outside~$O_1$.
Thus, the value of~$\varepsilon$ (as above) depends only on~$d_2$, which in turn depends only on~$d_1$.
Trajectories of~$X$ starting outside~$O_1$ therefore remain at distance at least~$\varepsilon$ from~$F$.
Hence $X$ is repelled by~$F$, as per Remark~\ref{rmk:repel_dist}, proving Claim~B.
\end{proof}

\subsection{Consequences, special cases, and clarif\/ications}
\label{sub:consequences}

Next, we give three corollaries of Theorem~\ref{t:vertexical}.
First, when the repulsing index set $R$ consists of all positive integers $\mathbb{Z}_{\geq1}$, the theorem specializes
to the following statement.
\begin{Corollary}
\label{c:repelled}
Fix a~vertexical family $\mathcal X=\{\mathcal X_S\}_{S\in\mathcal S}$ on open hypercubes indexed by the set~$\mathcal
S$ of all finite nonempty subsets of the positive integers $\mathbb{Z}_{\geq1}$.
If for every set~$S\in\mathcal S$, every differential inclusion $X\in\mathcal X_S$ is repelled by the vertex set
$\{0,1\}^S$ of its hypercube, then every such~$X$ is repelled by, and hence persistent relative to, the boundary
$\partial[0,1]^S$.
\end{Corollary}

In words, Corollary~\ref{c:repelled} states that to prove that every dif\/ferential inclusion in a~vertexical family is
repelled by the boundary, it suf\/f\/ices to show that each such dif\/ferential inclusion is repelled by the vertices.
The intuition behind this result is as follows.
If a~trajectory remains near a~proper face and away from its boundary for some time, then the vertexical property
allows us to project that part of the trajectory to a~trajectory of a~lower-dimensional dif\/ferential inclusion in
which the projected face is a~vertex, which is repelling by assumption; hence the original trajectory stays away from
the original face.

The next corollary is applied in our subsequent work~\cite{GeoGAC} to prove persistence results for dif\/ferential
inclusion families that arise from strongly endotactic reaction networks.
When the repulsing index set $R$ is empty, Theorem~\ref{t:vertexical}.1 specializes to the following statement.
\begin{Corollary}
\label{c:persistent}
Fix a~vertexical family $\mathcal X=\{\mathcal X_S\}_{S\in\mathcal S}$ on open hypercubes indexed by the set~$\mathcal
S$ of all finite nonempty subsets of the positive integers~$\mathbb{Z}_{\geq1}$.
Suppose that for every set $S\in\mathcal S$, every differential inclusion $X\in\mathcal X_S$ satisfies the
following hypotheses:
\begin{enumerate}\itemsep=-0.5pt
\item[$1)$] $X$ is repelled by the origin of its hypercube $(0,1)^{S}$, and
\item[$2)$] $X$ is persistent relative to the union of
all faces that do not contain the origin.
\end{enumerate}
Then every~$X\in\mathcal X_S$ is persistent relative to the boundary $\partial[0,1]^S$.
\end{Corollary}

Our f\/inal corollary is used to prove permanence-like results in our next work~\cite{GeoGAC}.
More precisely, Corollary~\ref{c:bounded} gives conditions under which trajectories of a~dif\/ferential inclusion~$X$
that begin in a~compact set~$K$ never leave a~larger compact set.
For ease of notation, we now introduce the following dif\/ferential inclusion.
\begin{Definition}
\label{d:restK}
In the setting of Def\/inition~\ref{d:persistent}, f\/ix a~dif\/ferential inclusion $X\subseteq TM$ and a~subset
$K\subseteq M$.
The \emph{restricted differential inclusion} $X_K\subseteq X$ is the smallest dif\/ferential inclusion such that
every trajectory of $X$ that begins in $K$ is a~trajectory of $X_K$.
\end{Definition}

In other words, $X_K$ consists of all tangent vectors to all trajectories of~$X$ that begin in~$K$.
\begin{Lemma}
\label{l:bounded}
Fix a~finite set $S$ and a~compact set $K\subseteq(0,1)^S$.
If a~differential inclusion $X\subseteq T(0,1)^S$ is repelled by the boundary~$\partial[0,1]^S$, then there exists
a~compact set $K^+$ with $K\subseteq K^+\subseteq(0,1)^S$ such that no trajectory of $X$ that begins in $K$ leaves~$K^+$.
\end{Lemma}
\begin{proof}
The open set $O_1=[0,1]^S\setminus K$ contains the boundary $\partial[0,1]^S$, so by def\/inition of repelling, there
exists an open set $O_2$ in $[0,1]^S$ that also contains the boundary $\partial[0,1]^S$ such that trajectories of $X_K$
that begin in $K$ never leave the compact set $K^+=[0,1]^S\setminus O_2$.
Hence, no trajectory of $X$ that begins in $K$ leaves $K^+$.
\end{proof}

\begin{Corollary}
\label{c:bounded}
Fix a~vertexical family $\mathcal X=\{\mathcal X_S\}_{S\in\mathcal S}$ on open hypercubes indexed by the set $\mathcal
S$ of all finite nonempty subsets of the positive integers $\mathbb{Z}_{\geq1}$, and a~repulsing index set
$R\subseteq\mathbb{Z}_{\geq1}$.
Assume the hypotheses of Theorem~{\rm \ref{t:vertexical}}.
Fix a~set $S\subseteq\mathcal S$, a~compact set $K\subseteq(0,1)^S$, and a~differential inclusion $X\in\mathcal X_S$.
If $X_K$ is repelled by the union of all opposite faces of\/~$[0,1]^S$, then there exists a~compact set $K^+$ with
$K\subseteq K^+\subseteq(0,1)^S$ such that no trajectory of $X$ that begins in $K$ leaves~$K^+$.
\end{Corollary}
\begin{proof}
Immediate from Theorem~\ref{t:vertexical}.2 and Lemma~\ref{l:bounded}.
\end{proof}

If there exists $d>0$ so that all trajectories of~$X$ starting in~$K$ remain at distance greater than~$d$ from all
opposite faces, then $X_K$ is repelled by the union of the opposite faces.
Consequently, in the context of Corollary~\ref{c:bounded}, it follows that there exists a~compact set $K^+$ such that
no trajectory of~$X$ starting in~$K$ leaves~$K^+$.
In particular, if the repulsing index set is the empty set, then the existence of such a~bound $d$ simply means an
upper bound $1-d\in(0,1)$ on all coordinate components of all trajectories of~$X$ beginning in~$K$.
\begin{Remark}
\label{r:promoted}
The signif\/icance of Corollary~\ref{c:bounded} is that, given the f\/low from~$K$, promoting its persistence relative
to the opposite faces to repulsion by the opposite faces results in its repulsion by the entire boundary.
In this form, it looks like a~weaker form of Corollary~\ref{c:repelled}, but restricted to those trajectories that
begin in~$K$.
Note that it is a~weaker form because in Corollary~\ref{c:bounded}, we assume not only that $X_K$ is repelled by all
vertices, but also that it is repelled by the union of all the opposite faces.
An assumption like this appears to be necessary: without this assumption, the projections of opposite faces are
vertices in lower-dimensional dif\/ferential inclusions that need not be repelling (only for $X_K$ are the opposite
faces assumed to be repelling).
\end{Remark}
\begin{Remark}
\label{r:vertexical}
In the statement of Theorem~\ref{t:vertexical}, if the goal is to prove that a~particular dif\/ferential inclusion
$X\in\mathcal X_S$ is persistent relative to~$\partial[0,1]^S$ (or repelled by the charged set), then it is enough to
assume a~slightly weaker hypothesis, namely that $(i)$~$X$ itself is persistent relative to the union of all opposite
faces and repelled by the charged vertices of the hypercube, and $(ii)$~the lower-dimensional set $\mathcal X_U$ for each
proper nonempty subset $U\subseteq S$ is persistent relative to the union of all opposite faces and repelled by the
charged vertices of its corresponding hypercube $[0,1]^U$.
\end{Remark}
\begin{Remark}
\label{r:why_can't_replace_repelling_with_persistent}
We devised the notion of ``repelled by'' (Def\/inition~\ref{d:persistent}) expressly for the purpose of proving
Theorem~\ref{t:vertexical}.
It is natural to ask whether this notion is necessary: is a~vertexical family that is persistent relative to the
charged vertices necessarily persistent relative to the boundary? In other words, in the statement of
Theorem~\ref{t:vertexical}, can instances of ``repelled by'' be replaced by ``persistent relative to''? The answer is
no: such a~replacement makes the theorem false.
Using the hypothesis that certain vertices are repelling, in the proof of Theorem~\ref{t:vertexical} one obtains
a~value $\varepsilon>0$ (for a~neighborhood of the face) that depends only on~$d_2/2$ (the thickness of the block of
the face).
However, under the weaker assumption of persistence relative to the vertices, this value~$\varepsilon$ depends on the
specif\/ic subinterval~$I'$.
A trajectory can repeatedly enter and exit the block, so that there are inf\/initely many relevant subintervals~$I'$.
In this case, it is possible for the trajectory to enter arbitrarily small neighborhoods without violating persistence.
Also recall from Remark~\ref{r:repelled==>persistent} that although repulsion implies persistence, the converse is
false.
\end{Remark}

\section{Reaction network theory}
\label{sec:CRNT}

In this section, we def\/ine reaction networks, reaction systems, and their properties.

\subsection{Reaction networks}
\label{sub:networks}

Our networks are more general than usual for chemical reaction network  theory~\cite{Fein79,HornJackson}.
\begin{Definition}
\label{def:crn}
Write $\mathit{OpnInt}=\big\{(a,b)\mid0\leq a~<b\leq\infty\big\}$ for the set of open subintervals of~$\mathbb{R}_{>0}$
and $\mathit{CmpctInt}=\big\{[a,b]\mid0<a~\leq b<\infty\big\}$ for the set of compact subintervals.
\begin{enumerate}\itemsep=0pt
\item A \emph{reaction network} $(S,\mathcal C,\mathcal R)$ is a~triple of f\/inite sets: a~set $S$ of \emph{species},
a~set $\mathcal C\subseteq\mathbb{R}^{S}$ of \emph{complexes}, and a~set $\mathcal R\subseteq\mathcal C\times\mathcal
C$ of \emph{reactions}.
\item The \emph{reaction graph} is the directed graph $(\mathcal C,\mathcal R)$ whose vertices are the complexes and
whose directed edges are the reactions.
\item A reaction $r=(y,y')\in\mathcal R$, also written $y\to y'$, has \emph{reactant}
$y=\operatorname{reactant}(r)\in\mathbb{R}^S$,
\emph{product} $y'=\target{r}\in\mathbb{R}^S$, and \emph{reaction vector}
\begin{gather*}
\flux{r}={\target{r}-\source{r}}.
\end{gather*}
\item The \emph{reaction diagram} is the realization $(\mathcal C,\mathcal R)\to\mathbb{R}^S$ of the reaction graph
that takes each reaction $r\in\mathcal R$ to the translate of $\flux r$ that joins $\source r$ to $\target r$.
\item A \emph{linkage class} is a~connected component of the reaction graph.
\end{enumerate}
\end{Definition}
\begin{Remark}
\label{rmk:genl_crn}
The chemical reaction network  theory literature usually imposes the following requirements for a~reaction network.
\begin{itemize}\itemsep=0pt
\item Each complex takes part in some reaction: for all $y\in\mathcal C$ there exists $y'\in\mathcal C$ such that
$(y,y')\in\mathcal R$ or $(y',y)\in\mathcal R$.
\item No reaction is trivial: $(y,y)\notin\mathcal R$ for all $y\in\mathcal C$.
\end{itemize}
Def\/inition~\ref{def:crn} does \emph{not} impose these conditions; in other words, our reaction graphs may include
isolated vertices or self-loops.
We drop these conditions to ensure that the projection of a~network~-- obtained by removing certain species~-- remains
a~network under our def\/inition even if some reactions become trivial (see Def\/inition~\ref{def:projective}.1).
In addition, like Craciun, Nazarov, and Pantea~\cite[\S~7]{CNP}, we allow arbitrary real complexes $y\in\mathbb{R}^S$,
so our setting is more general than that of usual chemical reaction networks, whose complexes $y\in\mathbb{Z}_{\geq0}^S$
are nonnegative integer combinations of species, as in the following def\/inition.
The ODE systems def\/ined in the next subsection that result from real complexes have been studied over the years and
called ``power-law systems'' (see Remark~\ref{r:usual}).
\end{Remark}

\begin{Definition}\label{def:network}
A reaction network $(S,\mathcal C,\mathcal R)$ is
\begin{enumerate}\itemsep=0pt
\item[1)] \emph{integer} if $\mathcal C\subseteq\mathbb{Z}^S$;
\item[2)] \emph{chemical} if $\mathcal C\subseteq\mathbb{Z}_{\geq0}^S$;
\item[3)] \emph{reversible} if the reaction graph of the network is undirected: a~reaction $(y,y')$ lies in~$\mathcal R$ if
and only if its reverse reaction $(y',y)$ also lies in~$\mathcal R$;
\item[4)] \emph{strongly connected} if the reaction graph of the network is strongly connected; that is, if the reaction
graph contains a~directed path between each pair of complexes;
\item[5)] \emph{weakly reversible} if every linkage class of the network is strongly connected.
\end{enumerate}
\end{Definition}
Note that a~network is strongly connected if and only if it is weakly reversible and has only one linkage class.

The next def\/initions introduce endotactic networks of~\cite{CNP}.
\begin{Definition}\label{def:partialOrder}
The standard basis of $\mathbb{R}^S$ indexed by~$S$ def\/ines a~canonical inner product $\langle\cdot,\cdot\rangle$ on~$\mathbb{R}^S$ with respect to which the standard basis is orthonormal.
Let $w\in\mathbb{R}^S$.
\begin{enumerate}\itemsep=0pt
\item The vector $w$ def\/ines a~preorder on $\mathbb{R}^S$, denoted by $\leq_w$, in which
\begin{gather*}
y\leq_w y'\;\iff\;\langle w,y\rangle\leq\langle w,y'\rangle.
\end{gather*}
Write $y<_w y'$ if $\langle w,y\rangle<\langle w,y'\rangle$.
\item For a~f\/inite subset $Y\subseteq\mathbb{R}^S$, denote by $\init wY$ the set of $\leq_w$-maximal elements of~$Y$:
\begin{gather*}
\init wY=\big\{y\in Y\,\big|\, \langle w,y\rangle\geq\langle w,y'\rangle~\text{for all}~y'\in Y\big\}.
\end{gather*}
\item For a~reaction network $(S,\mathcal C,\mathcal R)$, the set $\mathcal R_w\subseteq\mathcal R$ of
$w$-\emph{essential reactions} consists of those whose reaction vectors are not orthogonal to~$w$:
\begin{gather*}
\mathcal R_w=\big\{r\in\mathcal R\,\big|\,\langle w,\flux{r}\rangle\neq0\big\}.
\end{gather*}
\item The $w$-\emph{support} $\supp w{S,\mathcal C,\mathcal R}$ of the network is the set of vectors that are
$\leq_w$-maximal among reactants of $w$-essential reactions:
\begin{gather*}
\supp w{S,\mathcal C,\mathcal R}=\init w{\source{\mathcal R_w}}.
\end{gather*}
\end{enumerate}
\end{Definition}
\begin{Remark}
In order to simplify the computations in our next work~\cite{GeoGAC}, we dif\/fer from the usual
convention~\cite{CNP,Pantea}, by letting $\init wY$ denote the $\leq_w$-maximal elements rather than the
$\leq_w$-minimal elements.
Accordingly, the inequalities in Def\/inition~\ref{def:endotactic} are switched, so our def\/inition of endotactic is
equivalent to the usual one.
\end{Remark}
\begin{Definition}
\label{def:endotactic}
Fix a~reaction network $(S,\mathcal C,\mathcal R)$.
\begin{enumerate}\itemsep=0pt
\item The network $(S,\mathcal C,\mathcal R)$ is $w$-\emph{endotactic} for some $w\in\mathbb{R}^S$ if
\begin{gather*}
\langle w,\flux{r}\rangle<0
\end{gather*}
for all $w$-essential reactions $r\in\mathcal R_w$ such that $\source{r}\in\supp w{S,\mathcal C,\mathcal R}$.

\item The network $(S,\mathcal C,\mathcal R)$ is \emph{$W$-endotactic} for a~subset $W\subseteq\mathbb{R}^S$ if
$(S,\mathcal C,\mathcal R)$ is $w$-endotactic for all vectors $w\in W$.

\item The network $(S,\mathcal C,\mathcal R)$ is \emph{endotactic} if it is $\mathbb{R}^S$-endotactic.

\item\label{def:strong_endotactic}
$(S,\mathcal C,\mathcal R)$ is \emph{strongly endotactic} if it is endotactic and for every vector $w$ that is not
ortho\-go\-nal to the stoichiometric subspace of $(S,\mathcal C,\mathcal R)$, there exists a~reaction $y\to y'$ in
$\mathcal R$ such that
\begin{enumerate}\itemsep=0pt
\item[$(i)$] $y>_w y'$, and
\item[$(ii)$] $y$ is $\leq_w$-maximal among all reactants in $(S,\mathcal C,\mathcal R)$:
$y\in\init{w}{\source{\mathcal R}}$.
\end{enumerate}
\end{enumerate}

\end{Definition}
\begin{Remark}\label{r:endotactic}
Endotactic chemical reaction networks, which generalize weakly reversible networks, were introduced by Craciun,
Nazarov, and Pantea~\cite[\S~4]{CNP}.
Our def\/inition is slightly more general still, because we do not require the reaction networks to be chemical
(Def\/inition~\ref{def:endotactic}).
Strongly endotactic reaction networks are new; they give rise to strong results concerning persistence using our
techniques; see Theorem~\ref{thm:strong_endotactic}.
Strongly connected networks (i.e., weakly reversible networks with only one linkage class) are strongly endotactic.
\end{Remark}

\begin{Remark}
\label{r:endotactic_geometry}
For the geometric intuition behind Def\/inition~\ref{def:endotactic}, imagine a~hyperplane normal to $w$ that is
sweeping across the reaction diagram in~$\mathbb{R}^S$ from ``inf\/inity in direction~$w$''.
As this hyperplane sweeps, it stops when it f\/irst reaches the reactant $y$ of a~reaction $y\to y'$ that is not
perpendicular to~$w$.
If all such reactions do not point into the halfspace already swept by the hyperplane~-- that is, all such reactions
have product $y'$ outside of the open swept halfspace~-- then the network is $w$-endotactic.
Equivalently, the network is endotactic if no such reaction makes an acute angle with~$w$.
Illustrations can be found in~\cite{CNP}.

\looseness=-1
As for strongly endotactic networks, the sweeping hyperplane now stops when it f\/irst touches the reactant of any
reaction, whether or not it is perpendicular to~$w$.
Again we require that the products of all such reactions lie outside of the open swept halfspace, and in addition at
least one of these reactions is not perpendicular to~$w$.
If this condition is satisf\/ied for all vectors~$w$ not orthogonal to the stoichiometric subspace, then the network is
strongly endotactic.
Both endotactic and strongly endotactic networks capture the idea that extreme reactions should not point outward.
\end{Remark}
\begin{Example}
\label{e:endotactic}
Here we follow the usual convention of depicting a~network by its reaction graph or reaction diagram and writing
a~complex as, for example, $2A+B$ rather than $y=(2,1)$.
The Lotka--Volterra reaction network, consisting of the three reactions
\begin{gather*}
A\rightarrow2A,
\qquad
\qquad
A+B\rightarrow2B,
\qquad
\qquad
B\rightarrow0,
\end{gather*}
is not endotactic.
Reversing all three reactions yields the network
\begin{gather}
\label{eq:LV-rev}
2A\rightarrow A,
\qquad
\qquad
2B\rightarrow A+B,
\qquad
\qquad
0\rightarrow B,
\end{gather}
which is strongly endotactic, as can be verif\/ied from its reaction diagram:
$$
\includegraphics[scale=0.9]{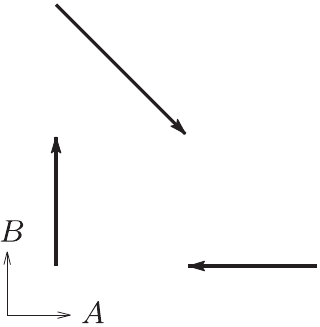}
$$
\end{Example}
\begin{Example}
Every weakly reversible reaction network is endotactic~\cite{CNP}.
However, even a~reversible reaction network may fail to be strongly endotactic, as in the following example of a~pair
of reversible reactions.
For $w=(0,-1)$, the $\leq_w$-maximal reactant complexes are at the bottom, but both of the corresponding reactions are
perpendicular to~$w$.
$$
\includegraphics[scale=0.9]{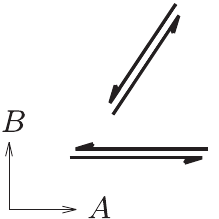}
$$
\end{Example}

\subsection{Reaction systems}
\label{sub:systems}
\begin{Definition}
\label{d:stoichiometry}
The \emph{stoichiometric subspace} $H$ of a~network is the span of its reaction vectors.
For a~positive vector $x_0\in\mathbb{R}_{>0}^S$, the \emph{invariant polyhedron} of $x_0$ is the polyhedron
\begin{gather}
\label{eq:scc}
\mathcal P=(x_0+H)\cap\mathbb{R}_{\geq 0}^S~.
\end{gather}
This polyhedron is also referred to as the \emph{stoichiometric compatibility class} in the chemical reaction network
theory literature~\cite{Fein79}.
\end{Definition}
\begin{Definition}
Let $(S,\mathcal C,\mathcal R)$ be a~reaction network.
\begin{enumerate}\itemsep=0pt
\item A \emph{tempering} is a~map $\kappa:\mathcal R\rightarrow\mathit{CmpctInt}$ that assigns to each reaction
a~nonempty compact positive interval.

\item A set $D\subseteq\mathbb{R}^S_{>0}$ is a~\emph{domain} if its intersection with every invariant
polyhedron~$\mathcal P$ of~$(S,\mathcal C,\mathcal R)$ is open in~$\mathcal P$.
\end{enumerate}
A \emph{reaction system} is a~triple consisting of a~reaction network, a~tempering, and a~domain.
\end{Definition}
\begin{Remark}\label{r:tempered}\looseness=-1
Mass-action dif\/ferential inclusions of reaction systems (Def\/inition~\ref{d:functor}) ge\-ne\-ra\-li\-ze the usual
mass-action kinetics ODE systems; see Remark~\ref{r:usual}.
One thinks of a~domain as a~promise that concentrations of species remain within the domain.
To explain the motivation behind temperings, recall that a~reaction network gives rise to a~dynamical system by way of
reaction rates.
For biochemical reaction networks, one is typically unable to measure precise values for the rates.
This occurs both because of incomplete information, and because of molecular and environmental variability.
One way to model this uncertainty is to allow reaction rates~$\kappa(r)$ to be time-dependent, as long as they are
uniformly bounded away from~$0$ and~$\infty$.
Craciun, Nazarov, and Pantea called such systems $\kappa$-variable~\cite{CNP}.
In a~similar spirit, we have chosen to work with dif\/ferential inclusions, allowing $\kappa(r)$ to take on every value
from an appropriate interval.
\end{Remark}
\begin{Definition}\qquad
\begin{enumerate}\itemsep=0pt
\item A \emph{confined} reaction system is a~reaction system whose domain is an invariant polyhedron of the
underlying reaction network.
\item An \emph{allotment} is a~map $\mu:S\to\mathit{OpnInt}$ sending each species $s\in S$ to an open positive interval.
The \emph{allotment hypercube} of an allotment $\mu$ is the open hypercube
$\square^\mu=\prod\limits_{s\in S}\mu(s)\subseteq\mathbb{R}_{>0}^S$.
A \emph{subconfined} reaction system is specif\/ied by a~reaction system and an allotment, in which the domain of the
reaction system is the intersection of the allotment hypercube and an invariant polyhedron of the underlying
reaction~network.
\end{enumerate}
\end{Definition}
\begin{Remark}
\label{r:confined_is_subconfined}
Every conf\/ined system is viewed as a~subconf\/ined system in which the allotment is understood to send every
species~$s$ to $(0,\infty)$, so the allotment hypercube is the entire positive orthant.
\end{Remark}
\begin{Remark}
The mathematical motivation prompting temperings and allotments is to ensure that projections of trajectories ``stay in
the family''.
Projections forget the exact concentrations of eliminated species.
Absorbing the ef\/fect of pre-projection dynamics into post-projection dynamics requires a~guarantee that the projected
species concentrations never leave a~certain suf\/f\/iciently large interval.
Therefore, in post-projection dynamics, the reaction rates remain within appropriately enlarged intervals.
Allotments and tempering reaction rates provide the extra f\/lexibility for this construction.

This intuition is made precise in Section~\ref{sec:funct}:\;interval-valued rates allow us to def\/ine an appropriate
domain category on which mass-action kinetics becomes functorial~({Theorem~\ref{thm:functor}}).
As a~consequence, projective classes of reaction systems (Def\/inition~\ref{def:projective}) give rise to families of
dif\/ferential inclusions that are vertexical.
\end{Remark}
\begin{Remark}
\label{r:confined}
Invariant polyhedra are forward-invariant sets with respect to the dynamics arising from mass-action kinetics; see
Remark~\ref{rmk:remain_in_scc}.
Therefore, a~conf\/ined reaction system allows us to restrict our attention to the dynamics on a~specif\/ic invariant
set.
\end{Remark}

\section{Functoriality of mass-action kinetics}
\label{sec:funct}

To every reaction system $N$ we assign a~dif\/ferential inclusion $\mathfrak{M}(N)$~(Def\/inition~\ref{d:functor}).
This assignment $\mathfrak{M}$ generalizes the usual mass-action kinetics ODE system in the chemical reaction network
theory literature~\cite{Fein79,HornJackson}.
It is the goal of this section to analyze how~$\mathfrak M$ behaves under projections of subconf\/ined reaction systems.
The main result (Theorem~\ref{t:proj-vertexical}) states that \emph{projective classes} of reaction
systems~(Def\/inition~\ref{def:projective}) give rise to vertexical families of dif\/ferential inclusions.
In particular, chemical, reversible, weakly reversible, endotactic, and strongly endotactic reaction systems all give
rise to vertexical families (Corollary~\ref{cor:isvertexical}).

The f\/irst task, which occupies Section~\ref{ss:cd}, is to make precise what is meant by projection, and by maps
between dif\/ferential inclusions.
It then becomes routine to verify two properties that are key to the proof of Theorem~\ref{t:proj-vertexical}, namely
that for every pair of subconf\/ined reaction systems $N_1$ and $N_2$ such that $N_2=p(N_1)$ is a~projection of~$N_1$,
the assignment $\mathfrak{M}$ induces a~map $\mathfrak{M}(p):\mathfrak{M}(N_1)\to\mathfrak{M}(N_1)$ between the
corresponding dif\/ferential inclusions such that
\begin{enumerate}\itemsep=0pt
\item[1)] the identity projection gets sent to the identity map on dif\/ferential inclusions, and
\item[2)] the composition
$p=p_2\circ p_1$ of two projection maps $p_1:N_1\to N_2$ and $p_2:N_2\to N_3$ gets sent to the composition
$\mathfrak{M}(p_2)\circ\mathfrak{M}(p_1)$ of the corresponding maps between dif\/ferential inclusions; that is,
\begin{gather*}
\mathfrak{M}(p)=\mathfrak{M}(p_2)\circ\mathfrak{M}(p_1).
\end{gather*}
\end{enumerate}

These two properties of $\mathfrak{M}$ are precisely the ones required by the def\/inition of a~\emph{functor} in
category theory.
Therefore, we f\/ind it economical to use this language~(Theorem~\ref{thm:functor}).
Readers unfamiliar with the language of category theory should read the word ``functor'' as shorthand for the two
properties.
This is the extent of the category theory used in~this~paper.

\subsection{Categorical def\/initions}
\label{ss:cd}
\begin{Definition}
\label{def:projective}
Recall, from Def\/inition~\ref{n:projectionmap}, the projection $\pi_U:\mathbb{R}^S\to\mathbb{R}^U$ for $U\subseteq S$,
and denote by $\pi_U^{\times2}=\pi_U\times\pi_U:\mathbb{R}^S\times\mathbb{R}^S\to\mathbb{R}^U\times\mathbb{R}^U$ the
product of $\pi_U$ with itself.
\begin{enumerate}\itemsep=0pt
\item For a~reaction network $(S,\mathcal C,\mathcal R)$, and a~nonempty subset $U\subseteq S$ of species, the
\emph{reduced reaction network} is the reaction network $\pi_U(S,\mathcal C,\mathcal R)=\big(U,\pi_U(\mathcal
C),\pi_U^{\times2}(\mathcal R)\big)$.
\item
\label{d:projective property}
A property $P$ of reaction networks is \emph{projective} if for all f\/inite nonempty sets $S$, reaction networks
$(S,\mathcal C,\mathcal R)$, and nonempty subsets $U\subseteq S$, if $(S,\mathcal C,\mathcal R)$ has property $P$ then
the reduced reaction network $\pi_U(S,\mathcal C,\mathcal R)$ has property~$P$.
\item The set of all reaction networks with a~given projective property is a~\emph{projective class}.
\end{enumerate}
\end{Definition}
\begin{Remark}
\label{r:reduced}
The reduced reaction network $\pi_U(S,\mathcal C,\mathcal R)$ is obtained from the reaction network $(S,\mathcal
C,\mathcal R)$ by deleting all species outside of~$U$.
This concept was def\/ined by Anderson~\cite[\S~3.2]{Anderson11}, who required that any trivial reactions be removed
from the reduced reaction set $(\pi_U\times\pi_U)(\mathcal R)$.
In contrast, we allow trivial reactions (Remark~\ref{rmk:genl_crn}).
Another related notion in the context of reversible reactions is that of ``reduced event-system'' introduced
in~\cite{adleman-2008}.
\end{Remark}
\begin{Example}
\label{e:reduced}
For $U=\{A\}$, the reduced network
\begin{gather*}
2A\to A\leftarrow0\leftarrow0
\end{gather*}
of the network~\eqref{eq:LV-rev} in Example~\ref{e:endotactic} is obtained by removing species~$B$.
The reduced network is strongly endotactic, as is the original network~\eqref{eq:LV-rev}.
\end{Example}

Next we see that the implication in Example~\ref{e:reduced} holds in general for strongly endotactic networks.
Such an implication was already completed for weakly reversible networks by Anderson~\cite[Lemma 3.4]{Anderson11}, and
for endotactic networks by Pantea~\cite[Proposition 3.1]{Pantea}.
\begin{Lemma}\label{lem:fams_proj_well}
The classes of integer, chemical, reversible, strongly connected, weakly reversible, endotactic, or strongly endotactic
reaction networks are projective.
Further, if~$P_1$ and~$P_2$ are projective properties, then so are the conjunction $P_1\wedge P_2$ and disjunction
$P_1\vee P_2$.
\end{Lemma}
\begin{proof}
Projectivity holds for integer and chemical networks because projection preserves integrality and nonnegativity of
points in~$\mathbb{R}^S$.
Projectivity holds for reversible, strongly connected, and weakly reversible networks because these conditions depend
only on the reaction graph, on whose vertices and edges projection is surjective.

Next, consider an endotactic network with species set~$S$ and the reduced network arising from a~nonempty subset
$U\subseteq S$.
Take any vector $w\in\mathbb{R}^U$.
For any reduced reaction $\pi_U(r)$, where $r$ is a~reaction in the original network,
\begin{gather}
\label{eq:fluxEq}
\big\langle w,\flux{\pi_U(r)}\big\rangle=\big\langle(w,0),\flux{r}\big\rangle.
\end{gather}
Thus, the $w$-essential reactions of the reduced network are the projections under $\pi_U^{\times2}$ of the
$(w,0)$-essential reactions of the original network, where we write $(w,0)\in\mathbb{R}^U\times\mathbb{R}^{S\setminus
U}$.
Similarly, the $w$-support of the reduced network is the projection under $\pi_U$ of the $(w,0)$-support of the
original network.
So, if $\pi_U(r)$ is a~$w$-essential reaction of the reduced network with reactant in the $w$-support, then the
original reaction $r$ is a~$(w,0)$-essential reaction of the original network with reactant in the $(w,0)$-support.
By~\eqref{eq:fluxEq} and the def\/inition of endotactic, $\big\langle
w,\flux{\pi_U(r)}\big\rangle=\big\langle(w,0),\flux{r}\big\rangle<0$.
Hence the reduced network is endotactic.

Next, let $H$ denote the stoichiometric subspace of a~strongly endotactic network, so $\pi_U(H)$ is the stoichiometric
subspace of the reduced network.
Take a~vector $w\in\mathbb{R}^{U}$ that is not orthogonal to $\pi_U(H)$.
We need only show that there exists a~reaction $\pi_U(y)\to\pi_U(y')$, where $y\to y'$ is a~reaction in the original
network, such that $\pi_U(y)>_w\pi_U(y')$ and $\pi_U(y)$ is $\leq_w$-maximal among all reactant vectors in the reduced
network.
Again consider $(w,0)\in\mathbb{R}^U\times\mathbb{R}^{S\setminus U}$.
As $w$ is not orthogonal to~$\pi_U(H)$, it follows that $(w,0)$ is not orthogonal to~$H$, and the preorder $\leq_w$ on
the reduced reactant complexes is the projection under $\pi_U$ of the preorder $\leq_{(w,0)}$ on the original reactant
complexes.
Since the original network is strongly endotactic, there is a~reaction $y\to y'$ in the original network with
$y>_{(w,0)}y'$ such that $y$ is $\leq_{(w,0)}$-maximal among all reactant vectors.
This reaction achieves our requirements.

The claim about conjunctions and disjunctions follows formally by Def\/inition~\ref{def:projective}.\ref{d:projective
property}.
\end{proof}
\begin{Notation}
\label{n:interval-arithmetic}
Let $I,J\subseteq\mathbb{R}_{\geq 0}$ be two intervals, and $S$ a~f\/inite nonempty set.
\begin{enumerate}\itemsep=0pt
\item Def\/ine $I\cdot J=\{i\cdot j\mid i\in I,j\in J\}\subseteq\mathbb{R}_{>0}$ and $\bigodot_{s\in S}I_s$
pointwise.

\item For $n\in\mathbb{Z}_{\geq1}$ def\/ine $I^n$ recursively as $I\cdot I^{n-1}$ with $I^1=I$.

\item $I\times J$ and $\prod_{s\in S}I_s$ denote Cartesian products of intervals, as usual.
\end{enumerate}
\end{Notation}
\begin{Definition}
\label{d:projectable}
Let $S$ be a~f\/inite nonempty set, and consider a~function $\mu:S\to\mathit{OpnInt}$ to the set of open positive
intervals.
A nonempty subset $U\subseteq S$ is \emph{$\mu$-projectable} if $0<\inf{\mu(s)}$ and $\sup{\mu(s)}<\infty$ for all
$s\in S\setminus U$; that is, the left and right endpoints of the intervals $\mu(s)$ are bounded away from $0$ and
$\infty$ for those $s$ outside of~$U$.
\end{Definition}
\begin{Remark}
\label{r:complement}
The condition of Def\/inition~\ref{d:projectable} is on the complement $S\setminus U$ because those are the species
removed in projecting to~$U$, and so it is those species that must be bounded away from~$0$ and~$\infty$.
The set $S$ itself is trivially $\mu$-projectable, for all $\mu:S\to\mathit{OpnInt}$.
\end{Remark}

We show that subconf\/ined reaction systems form a~category whose morphisms are projections, where the projection $p_U$
from one object~$N$ to another corresponds to substituting intervals from the allotment of~$N$ in place of a~set
$S\setminus U$ of forgotten species.
\begin{Definition}
\label{def:intWtd}
The category $\mathcal N$ of subconf\/ined reaction systems with projections is given by the following
data.
\begin{enumerate}\itemsep=0pt
\item Objects: each is a~subconf\/ined reaction system $N$, specif\/ied by a~reaction network $(S,\mathcal C,\mathcal
R)$ along with a~tempering $\kappa:\mathcal R\rightarrow\mathit{CmpctInt}$, an allotment $\mu:S\to\mathit{OpnInt}$, and
an invariant polyhedron  $\mathcal P=(x_0+H)\cap\mathbb{R}_{\geq 0}^{S}$.
 \item Morphisms: $p_U:N\to N'$
if
\begin{itemize}\itemsep=0pt
\item the network of~$N'$ is $(S',\mathcal C',\mathcal R')=\pi_U(S,\mathcal C,\mathcal R)$ for a~$\mu$-projectable
subset $U\subseteq S$; \item the tempering of~$N'$ is
$\displaystyle\kappa':\pi_U(r)\mapsto\kappa(r)\,\cdot\!\bigodot_{s\in S\setminus U}\mu(s)^{\source{r}_s}$, where the
exponent on $\mu(s)$ is the component indexed by~$s$ in the vector $\source r$; \item the allotment of~$N'$ is
$\mu'=\mu|_U$, gotten by restricting the allotment of~$N$ to~$U$;~and \item the invariant polyhedron  of~$N'$ is
$\mathcal P'=\left(\pi_U(x_0)+\pi_U(H)\right)\cap\mathbb{R}_{\geq 0}^{U}$.
\end{itemize}
\end{enumerate}
\end{Definition}
\begin{Remark}
\label{r:projected_scc}
In Def\/inition~\ref{def:intWtd}.2, $\pi_U(H)$ is the stoichiometric subspace of~$N'$ because it is spanned by the
vectors $\flux{\pi_U(r)}=\pi_U(\flux{r})$, where $r$ is a~reaction of~$N$.
Thus, $\left(\pi_U(x_0)+\pi_U(H)\right)\cap\mathbb{R}_{\geq 0}^{U}$ is an invariant polyhedron  of~$N'$.
\end{Remark}
\begin{Remark}
\label{r:composition}
Composition in $\mathcal N$ is well-def\/ined because f\/irst projecting to $U\subseteq S$ and then projecting to
$V\subseteq U$ is the same as projecting directly to~$V$.
\end{Remark}

Mass-action kinetics assigns to each subconf\/ined reaction system a~dif\/ferential inclusion on its domain.
Theorem~\ref{thm:functor} states that this assignment makes mass-action kinetics a~functor, with domain category
$\mathcal N$ and codomain category as follows.
\begin{Definition}
\label{d:DI}
The category $\mathcal{DI}$ of dif\/ferential inclusions is given by the following data.
\begin{enumerate}\itemsep=0pt
\item Objects: each is a~choice of manifold with corners and a~dif\/ferential inclusion on it.
\item Morphisms: a~morphism from $X\subseteq TM$ to $Y\subseteq TN$ is a~continuous map $k:M\to N$ such that for each
trajectory $f:I\to M$ of~$X$, there is a~trajectory $g:J\to N$ of~$Y$ and an order-preserving continuous map
$\alpha:I\to J$ satisfying
\begin{gather}
\label{eqn:1}
k\circ f=g\circ\alpha.
\end{gather}
\end{enumerate}
\end{Definition}
\begin{Lemma}\label{l:DI}
Composition of continuous maps induces a~well-defined composition on~$\mathcal{DI}$.
Specifically, assume $X_j\subseteq TM_j$ for $j=1,2,3$ are differential inclusions, with morphisms $k_{12}:M_1\to
M_2$ and $k_{23}:M_2\to M_3$ in~$\mathcal{DI}$.
If $f_1:I_1\to M_1$ is a~trajectory of~$X_1$, then there is a~trajectory $f_3:I_3\to M_3$ of~$X_3$ and an
order-preserving continuous map $\alpha_{13}:I_1\to I_3$ such that the composite continuous map $k_{13}=k_{23}\circ
k_{12}$ satisfies $k_{13}\circ f_1=f_3\circ\alpha_{13}$.
\end{Lemma}
\begin{proof}
Given $f_1$, since $k_{12}$ is a~morphism in~$\mathcal{DI}$, there is a~trajectory $f_2:I_2\to M_2$ of~$X_2$ and
a~continuous order-preserving map $\alpha_{12}:I_1\to I_2$ such that $k_{12}\circ f_1=f_2\circ\alpha_{12}$.
For the desired trajectory $f_3:I_3\to M_3$ of~$X_3$ use the one af\/forded by virtue of $k_{23}$ being a~morphism
in~$\mathcal{DI}$, given~$f_2$, which comes with a~continuous order-preserving map $\alpha_{23}:I_2\to I_3$ such that
$k_{23}\circ f_2=f_3\circ\alpha_{23}$.
Set $\alpha_{13}=\alpha_{23}\circ\alpha_{12}$.
Then
\begin{gather*}
k_{13}\circ f_1=k_{23}\circ k_{12}\circ f_1
=k_{23}\circ f_2\circ\alpha_{12}
=f_3\circ\alpha_{23}\circ\alpha_{12}
=f_3\circ\alpha_{13},
\end{gather*}
as desired.
\end{proof}
\begin{Remark}
\label{r:DI}
Compare the notion of morphism in~$\mathcal{DI}$ (Def\/inition~\ref{d:DI}) with that of vertexical family
(Def\/inition~\ref{d:vertexical}).
Equation~\eqref{eqn:1} also occurs in Def\/inition~\ref{d:vertexical}, with $k$ being a~particular type of continuous
map~$\pi_U$.
However, Def\/inition~\ref{d:DI} asks for a~global map~$k$, whereas the maps in Def\/inition~\ref{d:vertexical} are
required only locally, on blocks of faces.
\end{Remark}
\begin{Remark}
\label{r:topequiv}
One motivation for def\/ining the category of dif\/ferential inclusions this way comes from the dynamical systems
concept of \emph{topological equivalence}~\cite{irwin80smooth}, which identif\/ies two phase portraits as qualitatively
the same, even if the details of the dynamics may dif\/fer.
The isomorphisms in our category~$\mathcal{DI}$ correspond exactly to topological equivalence.

For this reason, morphisms between dif\/ferential inclusions may also be called \emph{topological morphisms}.
Intuitively, if a~topological morphism is a~monomorphism, then its target dif\/ferential inclusion qualitatively
simulates the domain dif\/ferential inclusion.
Maps that are not monomorphisms can of course result in the loss of information, in general.
The categorical message of Theorem~\ref{t:vertexical} is that it is sometimes possible to piece together many ``lossy''
maps on the same domain to regain substantial information about the domain~dynamics.
\end{Remark}

Another concept from dynamical systems, \emph{topological conjugacy}~\cite{irwin80smooth}, is a~stronger notion than
topological equivalence that disallows time reparameterization.
This motivates looking at a~subcategory $\mathcal{DI}_1$ of $\mathcal{DI}$ in which the order-preserving map $\alpha$
in Def\/inition~\ref{d:DI} is required to be the identity map; the def\/inition follows.
\begin{Definition}
\label{def:alphaIs1}
The \emph{category} $\mathcal{DI}_1$ \emph{of differential inclusions with topological semiconjugacy morphisms} is
the subcategory of~$\mathcal{DI}$ with the following data.
\begin{enumerate}\itemsep=0pt
\item Objects: the same objects as in $\mathcal{DI}$.
\item Morphisms: a~morphism from $X\subseteq TM$ to $Y\subseteq TN$ is a~continuous map $k:M\to N$ such that $k\circ f$
is a~trajectory of~$Y$ whenever $f:I\to M$ is a~trajectory of~$X$.
\end{enumerate}
\end{Definition}

The proof of Lemma~\ref{l:DI} makes it plain that any composition in~$\mathcal{DI}$ of morphisms in~$\mathcal{DI}_1$ is
a~morphism in~$\mathcal{DI}_1$, because every reparameterization map $\alpha_{ij}$ in that proof can be taken to be the
identity map on~$I_1$.

The next def\/inition uses the multinomial notation $x^y:=x_1^{y_1}x_2^{y_2}\cdots x_m^{y_m}$ for $x,y\in\mathbb{R}^m$.
\begin{Definition}
\label{d:functor}
The \emph{mass-action differential inclusion} of a~reaction system, specif\/ied by a~reaction network $(S,\mathcal
C,\mathcal R)$ with tempering~$\kappa$ and domain~$D$, is the dif\/ferential inclusion on $\mathbb{R}_{>0}^S$ whose
f\/iber over each point $x\in D$ is
\begin{gather}
\label{eq:m-a}
\bigg\{\sum_{r\in\mathcal R}k_r x^{\source{r}}\flux{r}\,\big|\, k_r\in\kappa(r)~\text{for all}~r\in\mathcal R\bigg\}\subseteq\mathbb{R}^S=T_x\mathbb{R}_{>0}^S
\end{gather}
and whose f\/iber over all points $x\in\mathbb{R}_{>0}^S\setminus D$ is empty.
The \emph{mass-action functor} $\mathfrak M:\mathcal N\to\mathcal{DI}_1$ from the category $\mathcal N$ of
subconf\/ined reaction systems to the category $\mathcal{DI}_1$ is def\/ined~on
\begin{enumerate}\itemsep=0pt
\item[1)] objects~$N\in\mathcal N$ by letting $\mathfrak M(N)$ be the mass-action dif\/ferential inclusion of~$N$, and~on
\item[2)] morphisms $p_U:N\to N'$ by letting $\mathfrak M(p_U)$ be the projection $\pi_U:\square^\mu\to\square^{\mu'}$.
\end{enumerate}
\end{Definition}
\begin{Remark}
\label{r:usual}
The usual mass-action kinetics ODE systems are the special cases in which the reaction network is chemical, the
tempering~$\kappa$ assigns not an interval of positive length but a~point~-- the reaction rate constant~-- to each
reaction, and the allotment hypercube is the entire positive orthant: $\square^\mu=\mathbb{R}_{>0}^S$.
The more general setting where stoichiometric coef\/f\/icients are allowed to be arbitrary (positive or negative) real
numbers has been studied in biochemical systems theory under the names ``power-law'' and ``generalized mass-action''.
This research area goes back to early work of Savageau in the 1960s~\cite{Savageau}.
\end{Remark}
\begin{Remark}
\label{rmk:m-a_in_cpc}
In the next subsection, we use a~dif\/feomorphism $\ell^S:\mathbb{R}_{>0}^S\to(0,1)^S$ to push forward mass-action
dif\/ferential inclusions to be def\/ined on hypercubes rather than orthants (Def\/inition~\ref{d:ell}).
Therefore, for mass-action dif\/ferential inclusions, the closure $\ol M$ of the allotment hypercube $M=\square^{\mu}$
is taken in the compactif\/ication $[0,\infty]^{S}$ when we are interested in the properties of being persistent,
permanent, or repelled.
\end{Remark}
\begin{Remark}
\label{rmk:remain_in_scc}
In a~mass-action dif\/ferential inclusion, only points in a~specif\/ied domain have nonempty f\/iber.
Imagine that we had instead def\/ined a~larger dif\/ferential inclusion, so that all f\/ibers in the positive orthant
are nonempty and take the form given by mass-action~\eqref{eq:m-a}.
Then, recalling that the stoichiometric subspace of a~network is spanned by the reaction vectors $\flux{r}$, it
follows from~\eqref{eq:m-a} that every trajectory of this larger dif\/ferential inclusion would be conf\/ined to some
invariant polyhedron~\eqref{eq:scc}.
(This also uses the fact that trajectories remain nonnegative~\cite[\S~2]{ADS11}.) That is, the invariant polyhedra are
forward-invariant with respect to the larger dif\/ferential inclusion.
Therefore it is in fact appropriate to restrict the dif\/ferential inclusion to a~given invariant polyhedron (as we did
in Def\/inition~\ref{d:functor}) and then to analyze properties such as persistence of the restriction.
\end{Remark}

\subsection{Functorial results and consequences}
\label{subsec:functor}
\begin{Theorem}
\label{thm:functor}
The mass-action functor $\mathfrak M$ is a~functor from the category $\mathcal N$ of subconfined reaction systems
with projection morphisms to the category $\mathcal{DI}_1$ of differential inclusions with topological semiconjugacy
morphisms.
\end{Theorem}
\begin{proof}
The content of the statement is twofold: f\/irst, that $\mathfrak M(p_U)=\pi_U$ in Def\/inition~\ref{d:functor}.2
indeed def\/ines a~morphism $\mathfrak M(N)\to\mathfrak M(N')$, and second, that $\mathfrak M$ preserves identity
morphisms as well as compositions.
The second is straightforward (projections are sent to projections).
For the f\/irst, it suf\/f\/ices to observe that for any nonempty $U\subseteq S$, the projection $\pi_U\circ f$ of
a~trajectory~$f$ of~$\mathfrak M(N)$ is a~trajectory of~$\mathfrak M(p_U(N))$ by def\/inition of~$p_U$.
This observation uses the fact that the image of~$f$ is in the invariant polyhedron~$\mathcal P$ of~$N$, so the image
of $\pi_U\circ f$ is in $\pi_U(\mathcal P)$ which is contained in the invariant polyhedron  of~$p_U(N)$.
\end{proof}
\begin{Remark}
\label{rem:allotvxical}
When we f\/ix an allotment, we do not, and can not, insist that trajectories stay within the allotment for all time.
Our assertion of functoriality is to the ef\/fect that for the period of time that a~trajectory does stay within the
allotment, its projection factors through a~smaller system.
The def\/inition of vertexical family~(Def\/inition~\ref{d:vertexical}) is weak enough to tolerate such a~weak
guarantee, and yet strong enough to be able to prove theorems on persistence and permanence.
In fact, the def\/inition of vertexical requires even less: the projection of a~trajectory is required to factor
through a~smaller system only for the period of time that it lies in a~corresponding block.
\end{Remark}

Vertexical families are def\/ined in terms of dif\/ferential inclusions on unit hypercubes $(0,1)^S$.
In contrast, the mass-action functor produces dif\/ferential inclusions on positive orthants~$\mathbb{R}_{>0}^S$.
To translate back and forth, our next def\/inition f\/ixes a~smooth, order-preserving dif\/feomorphism
$\mathbb{R}_{>0}\to(0,1)$.
The actual choice is irrelevant for our purposes, but if it helps, the reader may consider the function $x\mapsto
x/(1+x)$.
\begin{Definition}
\label{d:ell}
Fix a~smooth, order-preserving dif\/feomorphism $\ell:\mathbb{R}_{>0}\to(0,1)$.
For every nonempty f\/inite set~$S$, let $\ell^S:\mathbb{R}_{>0}^S\to(0,1)^S$, with derivative
$d\ell^S:T\mathbb{R}_{>0}^S\to T(0,1)^S$.
\begin{enumerate}\itemsep=0pt
\item The \emph{pushforward} of a~dif\/ferential inclusion $X$ on $\mathbb{R}_{>0}^S$ \emph{under~$\ell$} is the
dif\/ferential inclusion $d\ell^S(X)$ on $(0,1)^S$.

\item For a~subconf\/ined reaction system~$N$ def\/ined on a~species set~$S$ with allotment~$\mu$, the dif\/ferential
inclusion $\mathfrak M^\ell(N)$ is the pushforward of the mass-action dif\/ferential inclusion~$\mathfrak M(N)$,
considered as a~dif\/ferential inclusion on the image $\ell^S(\square^\mu)\subseteq(0,1)^S$ of the allotment hypercube
$\square^\mu$ under~$\ell$.
\end{enumerate}
\end{Definition}

The following is the main result of this section.
\begin{Theorem}
\label{t:proj-vertexical}
Fix a~projective property $P$ of reaction networks.
The class $\mathcal F_P$ of all confined reaction systems whose underlying reaction networks have property~$P$ yields
a~vertexical family $\mathfrak M^\ell(\mathcal F_P)=\big\{\mathfrak M^\ell(N)\mid N\in\mathcal F_P\big\}$ of
differential inclusions on open hypercubes.
\end{Theorem}
\begin{proof}
By def\/inition, $\mathfrak M^\ell(\mathcal F_P)$ is a~family of dif\/ferential inclusions on open hypercubes.
Hence we need only prove that $\mathfrak M^\ell(\mathcal F_P)$ is vertexical.
Consider a~{conf\/ined} reaction system $N\in\mathcal F_P$, specif\/ied by a~reaction network $(S,\mathcal C,\mathcal
R)$ with tempering~$\kappa$ and invariant polyhedron  $\mathcal P$.
Fix a~proper nonempty subset $U\subseteq S$ and a~vertex $x\in\{0,1\}^S$ of the hypercube $(0,1)^S$.
Denote by $F=F_{S\setminus U}(x)$ the corresponding face of the hypercube.

Fix $0<\varepsilon<\frac12$, and def\/ine $\mu'=\mu'_{\varepsilon,U}:S\to\mathit{OpnInt}$ by
\begin{gather*}
\mu'(s) =
\begin{cases} (0,\infty)&\text{if}~s\in U,
\\
\ell^{-1}(\varepsilon,1-\varepsilon)&\text{if}~s\in S\setminus U.
\end{cases}
\end{gather*}
Denote by~$N'$ the {subconf\/ined} reaction system $N'$ that agrees with~$N$ except that the allotment of~$N'$ is
$\mu'$.
The underlying reaction network of~$N'$ still has property~$P$, because $N$ and $N'$ have the same underlying reaction
network.

While $N'$ itself no longer has the entire positive orthant $\mathbb{R}_{>0}^S$ as its allotment hypercube, the set $U$
is $\mu'$-projectable (Def\/inition~\ref{d:projectable}), and the reduced network $p_U(N')$ has allotment
hypercube~$\mathbb{R}_{>0}^U$.
Since $P$ is projective, $p_U(N')$ therefore lies in~$\mathcal F_P$.
Consequently, by def\/inition of vertexical family (Def\/inition~\ref{d:vertexical}), it is enough to show that for any
trajectory $f:I\to F_\varepsilon$ of~$\mathfrak M^\ell(N)$ with image in the block~$F_\varepsilon$, the projection
$\pi_U\circ f=g$ is a~trajectory of the mass-action dif\/ferential inclusion $\mathfrak M^\ell(p_U(N'))$.

The trajectory $f$ is also a~trajectory of $\mathfrak M^\ell(N')$, since $F_\varepsilon\subseteq\ell^S(\square^{\mu'})$
by construction of~$\mu'$.
The result is now deduced easily from functoriality of mass-action kinetics (Theorem~\ref{thm:functor}): the morphism
$\pi_U$ in the category~$\mathcal{DI}_1$ from $\mathfrak M(N')$ to $\mathfrak M(p_U(N'))$ yields the pushforward
morphism $\mathfrak M^\ell(N')\to\mathfrak M^\ell(p_U(N'))$ in~$\mathcal{DI}_1$, so the proof is complete by
def\/inition of morphisms in the category $\mathcal{DI}_1$ (Def\/inition~\ref{def:alphaIs1}.2).
\end{proof}

Recall that in propositional logic, a~monotone (or monotonic) formula is one formed by the application of \textsc{and}
and \textsc{or} operations only, without the use of \textsc{not} operations.{\samepage
\begin{Corollary}
\label{cor:isvertexical}
Each of the classes of $($monotone combinations of$)$ integer, chemical, reversible, strongly connected, weakly reversible,
endotactic, or strongly endotactic confined reaction systems generates a~vertexical family of differential
inclusions.
\end{Corollary}
\begin{proof}
Immediate from Theorem~\ref{t:proj-vertexical} and Lemma~\ref{lem:fams_proj_well}.
\end{proof}}

\begin{Remark}
The category $\mathcal N$ of conf\/ined reaction systems with projections suf\/f\/ices for our purposes, but it would
be more natural to allow arbitrary reaction systems with an associated domain set $D\subseteq\mathbb{R}_{>0}^S$ that is
not necessarily derived from a~Cartesian product of intervals.
In addition, it is tempting to add to the category more morphisms, such as those corresponding to translation of the
reaction diagram within~$\mathbb{R}^S$, or scaling, rotation, arbitrary linear maps, graph homomorphisms, inversion
$z=1/x$, and so on.
It is easy to verify that translation acts as a~time-reparametrization.
Hence, even allowing translations, mass-action kinetics remains a~functor to $\mathcal{DI}$.
This leads to the following question, which we leave open.
\begin{Question}
\label{qopen}
What is the richest domain category for which mass-action kinetics remains a~functor to the category $\mathcal{DI}$ of
differential inclusions?
\end{Question}

This question is important because a~richer domain category would imply more ways of reducing the behavior of
a~network's mass-action kinetics to the behaviors of related networks.
This could allow us to ``program'' (and analyze) instances of reaction dynamics in high dimensions as appropriate
combinations of simpler reaction dynamics.
\end{Remark}
\begin{Remark}
Consider a~vertexical family of mass-action dif\/ferential inclusions for which the one-dimensional dif\/ferential
inclusions in the family are known to be permanent.
For instance, the dif\/ferential inclusions arising from weakly reversible, endotactic, or strongly endotactic networks
have this attribute.
It is tempting to attempt to argue that such a~family is permanent by the following induction: given a~trajectory, all
of its one-dimensional projections (which are also in the family, due to the functoriality of mass-action kinetics and
projectivity of the relevant properties) are permanent, and hence the trajectory itself eventually remains in a~compact
set.
However, this does not work because of uniformity issues.
It is true that for a~given trajectory, there exists a~compact set that it enters eventually.
However, what we need is one compact set so that every trajectory eventually enters this set, and our inductive
argument does not prove this.
A successful argument about permanence would require additional structure; for example, see Remark~\ref{r:lyapunov}.
\end{Remark}

\section{Implications for persistence of mass-action systems}
\label{sec:implications}

\looseness=-1
One of the long-standing open problems of chemical reaction network  theory is the global attractor conjecture
concerning so-called ``complex-balanced systems''.
Complex-balanced systems form a~well-studied subclass of weakly reversible mass-action ODE systems that contain all
so-called ``detailed-balanced'' systems and weakly reversible ``def\/iciency zero'' systems.
Many properties about complex-balanced systems (as well as detailed-balanced systems and def\/iciency zero systems)
were elucidated by Feinberg, Horn, and Jackson in the 1970s, and we provide only an overview here.
(See any of the references~\cite{TDS,Feinberg72,Fein79,HornJackson} for a~def\/inition of complex-balanced systems.)

For complex-balanced systems, it is known that there is a~unique steady state within the interior of each invariant
polyhedron  $\mathcal P$.
This steady state, called the Birch point in~\cite{TDS} due to the connection to Birch's theorem in algebraic
statistics, has a~strict Lyapunov function.
Therefore local asymptotic stability relative to $\mathcal P$ is guaranteed~\cite{Fein79, HornJackson} (see
Remark~\ref{r:lyapunov}).
An open question is whether all trajectories with an initial condition in the interior of $\mathcal P$ converge to the
unique Birch point of $\mathcal P$.
The assertion that the answer is ``yes'' is the content of the following conjecture, which was stated f\/irst by Horn
in 1974~\cite{horn74dynamics} and given the name ``Global Attractor Conjecture'' by Craciun et al.~\cite{TDS}.
\begin{Conjecture}[global attractor conjecture]
\label{c:GAC}
For any complex-balanced mass-action system and strictly positive initial condition~$x^0$, the Birch point in $\mathcal
P:=(x^0+H)\cap\mathbb{R}_{\ge0}^S$ $($see Definition~{\rm \ref{d:stoichiometry})} is a~global attractor of the relative interior of
the invariant polyhedron~$\operatorname{int}(\mathcal P)$.
\end{Conjecture}

In 1987, Feinberg~\cite{Fein87} conjectured the following.
\begin{Conjecture}[Feinberg's persistence conjecture]
\label{conj:persistence}
For every confined, weakly reversible mass-action ODE system $N$, the differential inclusion $\mathfrak M^\ell(N)$
is persistent.
\end{Conjecture}

Feinberg observed that the global attractor conjecture would follow from Conjecture~\ref{conj:persistence}, that is,
persistence of weakly reversible reaction networks.

Our functoriality results reduce persistence to bounding the dynamics away from vertices of the hypercube~-- so it
suf\/f\/ices to ensure that all species remain bounded away from~$0$ and~$\infty$ in the positive orthant~-- at the
price of considering tempered reaction systems.
More precisely, we state the following two corollaries.
\begin{Corollary}
\label{cor:repelled}
Let $\mathcal P$ be a~set of one or more of the following properties: integer, chemical, reversible, strongly
connected, weakly reversible, endotactic, and strongly endotactic.
Let $\mathcal F$ be the class of all confined reaction systems whose underlying reaction networks satisfy every
property in~$\mathcal P$.
If the mass-action differential inclusions of all reaction systems in~$\mathcal F$ are repelled by vertices after
pushing forward to open hypercubes by a~smooth order-preserving diffeomorphism, then they are repelled by the entire
boundary.
\end{Corollary}
\begin{proof}
Immediate from Corollaries~\ref{c:repelled} and~\ref{cor:isvertexical}.
\end{proof}

The next result states that another approach to persistence is by proving that trajectories are bounded and that the
origin is repelling; additionally, a~more uniform such bound yields a~permanence-like result.
\begin{Corollary}
\label{cor:persistent}
Let~$\mathcal P$ be a~set of one or more of the following properties: integer, chemical, reversible, strongly
connected, weakly reversible, endotactic, and strongly endotactic.
Let $\mathcal F$ be the class of all confined reaction systems whose underlying reaction networks satisfy every
property in~$\mathcal P$.
If the mass-action differential inclusions of all reaction systems in~$\mathcal F$ are repelled by the origin and
every trajectory of such a~differential inclusion is bounded, then
\begin{enumerate}\itemsep=0pt
\item[$1)$] these differential inclusions are persistent; and
\item[$2)$] if $X$ is such a~differential inclusion on
$\mathbb{R}_{>0}^S$, and $K\subseteq\mathbb{R}_{>0}^S$ is a~compact set for which there exists $A\in\mathbb{R}_{>0}$ such
that every trajectory of~$X$ that starts in~$K$ remains bounded above by~$A$ in each coordinate for all time, then for
some compact set $K^+\subseteq\mathbb{R}_{>0}^S$, no trajectory of~$X$ that begins in~$K$ leaves~$K^+$.
\end{enumerate}
\end{Corollary}
\begin{proof}
A dif\/ferential inclusion on $\mathbb{R}_{>0}^{S}$ is repelled by the origin if and only if after pushing forward to
open hypercubes the resulting dif\/ferential inclusion is repelled by the origin: the open sets in the def\/inition of
repelled move between hypercubes and positive orthants via the dif\/feomorphism.
Additionally, persistence of a~mass-action dif\/ferential inclusion is viewed with respect to the compactif\/ication
$[0,\infty]^{S}$, so persistence is equivalent to persistence of the pushforward with respect to $[0,1]^{S}$.
Finally, compact sets $K$ and bounds $A$ also move between hypercubes and positive orthants via the dif\/feomorphism.
Thus, the conclusion follows from Corollaries~\ref{cor:isvertexical},~\ref{c:persistent}, and~\ref{c:bounded} (for
which the repelled set is taken to be the origin~-- see also the description after Corollary~\ref{c:bounded}).
\end{proof}

Recently, Craciun, Nazarov, and Pantea~\cite{CNP} generalized Feinberg's persistence conjecture
(Conjecture~\ref{conj:persistence}) in the following three ways: the weakly reversible hypothesis is weakened to
endotactic, f\/ixed reaction rate constants are allowed to vary within bounded intervals (i.e., a~tempering), and the
conclusion of persistence is strengthened to permanence~\cite[\S~4]{CNP}.
\begin{Conjecture}[extended permanence conjecture]
\label{conj:perm}
For every confined endotactic reaction system $N$, the differential inclusion $\mathfrak M^\ell(N)$ is permanent.
\end{Conjecture}
\begin{Remark}
\label{r:projective_class}
Endotactic networks constitute a~projective class (Corollary~\ref{cor:isvertexical}), and, in our view, it is this
property that allowed Craciun, Nazarov, and Pantea to make their projection-type arguments~\cite{CNP, Pantea}.
Similarly, the property of having only one linkage class is projective: a~network with only one linkage class maintains
this property after reduction.
This type of argument was used in Anderson's proof of the global attractor conjecture for networks possessing only one
linkage class~\cite{Anderson11}.
Indeed, our work was motivated in part by the works both of Anderson and of Craciun, Nazarov, and Pantea.
\end{Remark}
\begin{Remark}
\label{r:lyapunov}
Horn and Jackson~\cite{HornJackson} established that any complex-balanced mass-action ODE system~$N$ on a~set~$S$ of
species admits the strict Lyapunov function
\begin{gather*}
g_\alpha(x)=\sum_{i\in S}x_i\left(\log\frac{x_i}{\alpha_i}-1\right),
\end{gather*}
where $\alpha\in\mathbb{R}_{>0}^S$ is a~given complex-balanced steady state of $N$.
Consequently, the dif\/ferential inclusion $\mathfrak M^\ell(N)$ is repelled by the vertex $0$, and its trajectories
are bounded away from faces not incident to $0$.
The proof of this assertion, which is due in part to Anderson and Craciun et al.~\cite{Anderson08,TDS}, proceeds by
showing that the Lyapunov function $g_\alpha$ has a~local maximum at the vertex~$0$ as well as bounded level sets that
do not allow trajectories to escape to inf\/inity.
If the class of complex-balanced systems were projective, then our arguments would have proved the global attractor
conjecture.
However, this is not the case: projections of complex-balanced systems need not be complex-balanced.
Therefore, Theorem~\ref{t:proj-vertexical} does not apply: we can not prove that complex-balanced systems form
a~vertexical family.
In fact, it can be shown that this family is not vertexical.
On the other hand, weakly reversible or endotactic networks do def\/ine a~vertexical family; recall
Remark~\ref{r:projective_class}.
\end{Remark}

The motivation for our work was to make progress on Conjectures~\ref{c:GAC},~\ref{conj:persistence},
and~\ref{conj:perm}.
The next two theorems state what we have accomplished in this direction.
First, we show that an extension of Feinberg's persistence conjecture in dimension $n$ implies the global attractor
conjecture in dimension $n+1$, under an additional assumption that the origin is repelling.
\begin{Theorem}
\label{t:n=>n+1}
Let $n$ be a~positive integer.
If for every confined, weakly reversible reaction system $N$ with no more than $n$ species, both
\begin{enumerate}\itemsep=0pt
\item[$1)$] $\mathfrak M^\ell(N)$ is persistent, and
\item[$2)$] $\mathfrak M^\ell(N)$ is repelled by the origin,
\end{enumerate}
then the global attractor conjecture $($Conjecture~{\rm \ref{c:GAC})} holds for complex-balanced mass-action systems with $n+1$
or fewer species.
\end{Theorem}
\begin{proof}\looseness=1
Let $\mathcal X$ denote the vertexical family of all dif\/ferential inclusions that arise from endotactic networks.
Let $N$ denote a~conf\/ined, complex-balanced mass-action system with $|S|\leq n+1$.
As mentioned after Conjecture~\ref{conj:persistence}, it suf\/f\/ices to prove that $N$ is persistent.
By Remark~\ref{r:vertexical} in the context of the empty repulsing index set, to prove that $N$ is persistent, it
suf\/f\/ices to show that $(i)$~$\mathfrak M^\ell(N)$ itself is persistent relative to all faces of $[0,1]^S$ that do not
meet the origin and repelled by the origin of the hypercube, and $(ii)$~the lower-dimensional set $\mathcal X_U$ for
each proper nonempty subset $U\subseteq S$ is persistent relative to all faces of $[0,1]^U$ that do not meet the origin
and repelled by the origin of the hypercube $[0,1]^U$.
Claim~$(i)$ follows from the Lyapunov function, as explained in Remark~\ref{r:lyapunov}: persistence relative to all
faces that do not meet the origin is equivalent to having bounded trajectories.
As for claim~$(ii)$, endotactic systems having lower dimension than~$X$ are persistent by hypothesis, and repelled by the
origin also by~hypothesis.
\end{proof}

The next theorem reduces the global attractor conjecture and Feinberg's persistence conjecture to the following
assertion that \emph{for mass-action differential inclusions arising from weakly reversible networks, trajectories
are bounded and repelled by the origin}.
\begin{Conjecture}
\label{conj:wr_are_repelled_bounded}
For every confined, weakly reversible reaction system $N$, the differential inclusion $\mathfrak{M}(N)$ is repelled
by the origin, and every trajectory of~$\mathfrak{M}(N)$ is bounded.
\end{Conjecture}
\begin{Theorem}
\label{t:bdd_tra_orig_rep=>conjs}
Conjecture~{\rm \ref{conj:wr_are_repelled_bounded}} implies the global attractor conjecture $($Conjecture~{\rm \ref{c:GAC})}.
\end{Theorem}
\begin{proof}
Consider a~conf\/ined weakly reversible reaction system $N$.
By hypothesis, the dif\/ferential inclusion $\mathfrak{M}(N)$ is repelled by the origin, and every trajectory of
$\mathfrak{M}(N)$ is bounded, so $\mathfrak{M}(N)$ is persistent by Corollary~\ref{cor:persistent}.1.
Thus, all weakly reversible systems are persistent, which implies the global attractor conjecture: see the discussion
after Conjecture~\ref{conj:persistence}.
\end{proof}

Conjectures~\ref{c:GAC},~\ref{conj:persistence},~\ref{conj:perm}, and~\ref{conj:wr_are_repelled_bounded} all remain
open; for an overview of recent progress on these problems, we refer the reader to the work of Anderson~\cite[\S~1.1]{Anderson08}.
However, using the results presented in the current paper, we prove the following in a~subsequent paper~\cite{GeoGAC},
which extends recent results of Anderson~\cite{Anderson11}.
\begin{Theorem}
\label{thm:strong_endotactic}
For every confined strongly endotactic reaction system $N$, the differential inclusion $\mathfrak M^\ell(N)$ is
permanent.
\end{Theorem}
\begin{proof}
To be proved in a~subsequent paper~\cite{GeoGAC}.
\end{proof}

The key step contributed by the results in the current paper is Corollary~\ref{cor:persistent} in the case of strongly
endotactic networks.
The approach in~\cite{GeoGAC} shows that for the mass-action dif\/ferential inclusions of strongly endotactic networks,
outside a~compact set the function $g_\alpha(x)$ in Remark~\ref{r:lyapunov} continues to behave like a~Lyapunov
function.
Theorem~\ref{thm:strong_endotactic} is established by an argument along the lines suggested in Remark~\ref{r:lyapunov}.

\subsection*{Acknowledgements} MG was supported by a~Ramanujan fellowship from the Department of Science and
Technology, India, and, during a~semester-long stay at Duke University, by the Duke MathBio RTG grant NSF DMS-0943760.
EM had support from NSF grant DMS-1001437.
AS was supported by an NSF postdoctoral fellowship DMS-1004380.
The authors thank David~F.~Anderson, Gheorghe Craciun, and Casian Pantea for helpful discussions, and Duke University
where many of the conversations occurred.
The authors also thank the two referees, whose perceptive and insightful comments improved this work.

\pdfbookmark[1]{References}{ref}
\LastPageEnding


\begin{thebibliography}{99}
\footnotesize\itemsep=0pt

\bibitem{adleman-2008}
Adleman L., Gopalkrishnan M., Huang M.D., Moisset P., Reishus D., On the
  mathematics of the law of mass action, \href{http://arxiv.org/abs/0810.1108}{arXiv:0810.1108}.

\bibitem{Anderson11}
Anderson D.F., A proof of the global attractor conjecture in the single linkage
  class case, \href{http://dx.doi.org/10.1137/11082631X}{\textit{SIAM~J. Appl. Math.}} \textbf{71} (2011), 1487--1508,
  \href{http://arxiv.org/abs/1101.0761}{arXiv:1101.0761}.

\bibitem{Anderson08}
Anderson D.F., Global asymptotic stability for a class of nonlinear chemical
  equations, \href{http://dx.doi.org/10.1137/070698282}{\textit{SIAM~J. Appl. Math.}} \textbf{68} (2008), 1464--1476,
  \href{http://arxiv.org/abs/0708.0319}{arXiv:0708.0319}.

\bibitem{AndersonShiu10}
Anderson D.F., Shiu A., The dynamics of weakly reversible population processes
  near facets, \href{http://dx.doi.org/10.1137/090764098}{\textit{SIAM~J. Appl. Math.}} \textbf{70} (2010), 1840--1858,
  \href{http://arxiv.org/abs/0903.0901}{arXiv:0903.0901}.

\bibitem{angeli2007Petri}
Angeli D., De~Leenheer P., Sontag E., A {P}etri net approach to persistence
  analysis in chemical reaction networks, in Biology and Control Theory:
  Current Challenges, \href{http://dx.doi.org/10.1007/978-3-540-71988-5_9}{\textit{Lect. Notes Contr. Inf.}}, Vol.~357, Editors
  I.~Queinnec, S.~Tarbouriech, G.~Garcia, S.I.~Niculescu, Springer-Verlag,
  Berlin, 2007, 181--216, \href{http://arxiv.org/abs/q-bio/0608019}{q-bio/0608019}.

\bibitem{ADS11}
Angeli D., De~Leenheer P., Sontag E.D., Persistence results for chemical
  reaction networks with time-dependent kinetics and no global conservation
  laws, \href{http://dx.doi.org/10.1137/090779401}{\textit{SIAM~J. Appl. Math.}} \textbf{71} (2011), 128--146.

\bibitem{Aubin}
Aubin J.P., Cellina A., Dif\/ferential inclusions. Set-valued maps and viability
  theory, \href{http://dx.doi.org/10.1007/978-3-642-69512-4}{\textit{Grundlehren der Ma\-the\-matischen Wissenschaften}}, Vol.~264,
  Springer-Verlag, Berlin, 1984.

\bibitem{BanajiMier}
Banaji M., Mierczy{\'n}ski J., Global convergence in systems of dif\/ferential
  equations arising from chemical reaction networks, \href{http://dx.doi.org/10.1016/j.jde.2012.10.018}{\textit{J.~Differential
  Equations}} \textbf{254} (2013), 1359--1374, \href{http://arxiv.org/abs/1205.1716}{arXiv:1205.1716}.

\bibitem{TDS}
Craciun G., Dickenstein A., Shiu A., Sturmfels B., Toric dynamical systems,
  \href{http://dx.doi.org/10.1016/j.jsc.2008.08.006}{\textit{J.~Symbolic Comput.}} \textbf{44} (2009), 1551--1565,
  \href{http://arxiv.org/abs/0708.3431}{arXiv:0708.3431}.

\bibitem{CNP}
Craciun G., Nazarov F., Pantea C., Persistence and permanence of mass-action
  and power-law dynamical systems, \href{http://dx.doi.org/10.1137/100812355}{\textit{SIAM~J. Appl. Math.}} \textbf{73}
  (2013), 305--329, \href{http://arxiv.org/abs/1010.3050}{arXiv:1010.3050}.

\bibitem{Fein87}
Feinberg M., Chemical reaction network structure and the stability of complex
  isothermal reactors: I.~The def\/iciency zero and def\/iciency one theorems,
  \href{http://dx.doi.org/10.1016/0009-2509(87)80099-4}{\textit{Chem. Eng. Sci.}} \textbf{42} (1987), 2229--2268.

\bibitem{Feinberg72}
Feinberg M., Complex balancing in general kinetic systems, \href{http://dx.doi.org/10.1007/BF00255665}{\textit{Arch.
  Rational Mech. Anal.}} \textbf{49} (1972), 187--194.

\bibitem{Fein79}
Feinberg M., Lectures on chemical reaction networks, unpublished lecture notes,
  1979, available at
 \mbox{\url{http://www.che.eng.ohio-state.edu/~FEINBERG/LecturesOnReactionNetworks/}}.

\bibitem{Freed}
Freedman H.I., Moson P., Persistence def\/initions and their connections,
  \href{http://dx.doi.org/10.2307/2048133}{\textit{Proc. Amer. Math. Soc.}} \textbf{109} (1990), 1025--1033.

\bibitem{GeoGAC}
Gopalkrishnan M., Miller E., Shiu A., A geometric approach to the global
  attractor conjecture, in preparation.

\bibitem{horn74dynamics}
Horn F., The dynamics of open reaction systems, in Mathematical Aspects of
  Chemical and Biochemical Problems and Quantum Chemistry ({P}roc. {SIAM}-{AMS}
  {S}ympos. {A}ppl. {M}ath., {N}ew {Y}ork, 1974), \textit{SIAM-AMS
  Proceedings}, Vol.~8, Amer. Math. Soc., Providence, R.I., 1974, 125--137.

\bibitem{HornJackson}
Horn F., Jackson R., General mass action kinetics, \href{http://dx.doi.org/10.1007/BF00251225}{\textit{Arch. Rational Mech.
  Anal.}} \textbf{47} (1972), 81--116.

\bibitem{irwin80smooth}
Irwin M.C., Smooth dynamical systems, \textit{Pure and Applied Mathematics},
  Vol.~94, Academic Press Inc., New York, 1980.

\bibitem{Lee02smooth}
Lee J.M., Introduction to smooth manifolds, \textit{Graduate Texts in
  Mathematics}, Vol.~218, Springer-Verlag, New York, 2003.

\bibitem{Pantea}
Pantea C., On the persistence and global stability of mass-action systems,
  \href{http://dx.doi.org/10.1137/110840509}{\textit{SIAM~J. Math. Anal.}} \textbf{44} (2012), 1636--1673,
  \href{http://arxiv.org/abs/1103.0603}{arXiv:1103.0603}.

\bibitem{Savageau}
Savageau M.A., Biochemical systems analysis: I.~Some mathematical properties of
  the rate law for the component enzymatic reactions, \href{http://dx.doi.org/10.1016/S0022-5193(69)80026-3}{\textit{J.~Theor. Biol.}}
  \textbf{25} (1969), 365--369.

\bibitem{Siegel}
Siegel D., Johnston M.D., A stratum approach to global stability of complex
  balanced systems, \href{http://dx.doi.org/10.1080/14689367.2010.545812}{\textit{Dyn. Syst.}} \textbf{26} (2011), 125--146,
  \href{http://arxiv.org/abs/1008.1622}{arXiv:1008.1622}.

\end{thebibliography}
\end{document}